\documentclass[12pt]{article}

\usepackage{amsfonts}
\usepackage{color}
\usepackage{graphicx}
\usepackage{amssymb}
\usepackage{listings}
\usepackage{float}
\usepackage{amsmath}
\usepackage{indentfirst}

\def\sqr#1#2{{\vcenter{\vbox{\hrule height.#2pt
				\hbox{\vrule width.#2pt height#1pt \kern#1pt \vrule width.#2pt}
				\hrule height.#2pt}}}}
\def\signed #1{{\unskip\nobreak\hfil\penalty50
		\hskip2em\hbox{}\nobreak\hfil#1
		\parfillskip=0pt \finalhyphendemerits=0 \par}}
\def\endpf{\signed {$\sqr69$}}

\def\dbR{{\mathop{\rm l\negthinspace R}}}

\def\3n{\negthinspace \negthinspace \negthinspace }
\def\2n{\negthinspace \negthinspace }
\def\1n{\negthinspace }

\def\ds{\displaystyle}

\def\dbR{{\mathop{\rm l\negthinspace R}}}
\def\={\buildrel \triangle \over =}

\def\a{\alpha}
\def\b{\beta}
\def\g{\gamma}

\def\e{\varepsilon}

\def\k{\kappa}
\def\l{\lambda}

\def\n{\nabla}

\def\t{\times}

\def\th{\theta}
\def\o{\omega}

\def\i{\infty}
\def\ns{\noalign{\ss} }

\def\G{\Gamma}
\def\D{\Delta}
\def\Th{\Theta}

\def\O{\Omega}

\def\cP{{\cal P}}

\def\cl{{\cal l}}
\def\no{\noindent}

\def\ms{\medskip}

\def\q{\quad}
\def\qq{\qquad}

\def\max{\mathop{\rm max}}
\def\min{\mathop{\rm min}}

\def\pa{\partial}
\def\h{\widehat}
\def\wt{\widetilde}
\def\cd{\cdot}

\def\supp{\hbox{\rm supp$\,$}}

\def\cl{\overline}

\def\Re{{\mathop{\rm Re}\,}}

\def\|{\Big |}
\def\({\Big (}
\def\){\Big )}
\def\[{\Big[}
\def\]{\Big]}
\def\be{\begin{equation*}}
\def\bel{\begin{equation}\label}
\def\ee{\end{equation}}
\def\bt{\begin{theorem}}
\def\bcd{\begin{condition}}
\def\ecd{\end{condition}}
\def\et{\end{theorem}}
\def\bc{\begin{corollary}}
\def\ec{\end{corollary}}
\def\bde{\begin{definition}}
\def\ede{\end{definition}}
\def\bl{\begin{lemma}}
\def\el{\end{lemma}}
\def\bp{\begin{proposition}}
\def\ep{\end{proposition}}
\def\br{\begin{remark}}
\def\er{\end{remark}}
\def\ba{\begin{array}}
\def\ea{\end{array}}
\def\ed{\end{document}}
\def\ns{\noalign{\ms}}
\def\ds{\displaystyle}

\def\square#1{\vbox{\hrule\hbox{\vrule height#1%
			\kern#1\vrule}\hrule}}
\def\rectangle#1#2{\vbox{\hrule\hbox{\vrule height#1%
			\kern#2\vrule}\hrule}}

\font\tenbb=msbm10 \font\sevenbb=msbm7 \font\fivebb=msbm5

\newfam\bbfam
\scriptscriptfont\bbfam=\fivebb \textfont\bbfam=\tenbb
\scriptfont\bbfam=\sevenbb

\def\oO{{\overline \O}}

\newtheorem{lemma}{Lemma}[section]
\newtheorem{remark}{Remark}[section]

\newtheorem{theorem}{Theorem}[section]
\newtheorem{corollary}{Corollary}[section]

\newtheorem{definition}{Definition}[section]
\newtheorem{proposition}{Proposition}[section]
\newtheorem{condition}{Condition}[section]

\makeatletter

\@addtoreset{equation}{section}
\makeatother
\begin{document}
\title{\bf Matlab program method of computing Carleman estimates and applications\thanks{This work is partially supported by the NSF of China
		under grant 11971333, 11931011, and by the Science Development Project of Sichuan University under grant 2020SCUNL201. }}
		
\author{Xiaoyu Fu\thanks{School
		of Mathematics, Sichuan University, Chengdu
		610064, China.  xiaoyufu@scu.edu.cn}, \  Yuan Gao\thanks{School
		of Mathematics, Sichuan University, Chengdu
		610064, China.  mathgaoyuan@stu.scu.edu.cn} \
 and \ Qingmei Zhao\thanks{  Corresponding author. School of Mathematics, Sichuan University, Chengdu 610064, China. E-mail address: zhaoqm@stu.scu.edu.cn}}

\date{}

\maketitle

\begin{abstract}

In this paper, we introduce a Matlab program method to compute Carleman estimate for the fourth order partial differential operator $\g\pa_t+\pa_x^4\ (\g\in\dbR)$. We obtain two kinds of Carleman estimates with different weight functions, i.e.  singular weight function and regular weight function, respectively. Based on  Carleman estimate with singular weight function,  one can obtain the known controllability and observability results for the 1-d fourth order parabolic-type equation, while based on Carleman estimate with regular weight function, one can deduce not only the known result on conditional stability in the inverse problem of half-order fractional diffusion equation, but also a new result on conditional stability in the inverse problem of half-order  fractional Schr\"odinger equation.

\end{abstract}

\no{\bf Key Words}. fourth order  parabolic-type operator, Carleman estimate, Matlab program method,  controllability, inverse problem.

\section{Introduction and main results }\label{s1}

In 1939,  Carleman \cite{c} established the strong unique continuation property for second order elliptic operators in dimensional two. The powerful technique he used, the so called ``Carleman weighted inequality"
has become one of the major tools in the study of unique continuation property, control problems and inverse problems for partial differential equations (see \cite{FLZ, FI, Hormander, Isakov2, Kn, Kli,  Y} and the rich references therein).

Roughly speaking, Carleman estimate can be regarded as weighted energy inequality. We give an example to explain this. Let $\Omega$ be a connected open set in $\dbR^{n}$, and let $P=P(x,D)$ be a differential operator of order $m$ in $\Omega$. Assume that there is a suitable function $\phi(\cd)\in C^\i(\oO;\dbR)$ satisfying $\ds \n\phi(x)\neq 0$, $x\in\O$, we say that the Carleman type estimate holds for $P$ if there exists a constant $C>0$ such that
  $$\ba{ll}\ds
  \sum_{|\alpha|<m}\l^{2(m-|\alpha|)-1}\int_{\Omega}|D^{\alpha}u|^{2}e^{2\l\phi}dx\leq C\int_{\Omega}|P(x,D)u|^{2}e^{2\l\phi}dx,
  \ea
  $$
where $u\in C_0^\i(\O)$, $\l$ is a parameter which can be chosen large enough.

Up to now, there are numerous results on Carleman estimates for the second order partial differential operators, the corresponding applications are well understood (for example, see \cite{FLZ} and the references cited therein). In the case of Carleman estimates for the higher order partial differential operators, in this respect, we refer to  \cite{CM, XCY, ZZ}  for the fourth order  parabolic-type operators in $1-d$ case, \cite{GK} for the operator $-\pa_t+\D^2$ (in $n-d$ case), \cite{HK} for the operator $\pa_t-A^2$ with $A=\n\cdot (a(x)\n )$ (in $n-d$ case), and \cite{Zheng} for the fourth order Schr\"odinger operator $i\pa_t+\pa_x^4$. Compared to Carleman estimate for second order partial differential operator, the computation of Carleman estimate for higher order partial differential operator is not an easy task, especially for the higher order partial differential operator in multi-dimensional case.

To simplify  complex computations, we introduce a Matlab program method for the proof of  Carleman estimate, via which one can obtain the desired result quickly and accurately. To see the idea of this method clearly, in this paper we  consider Carleman estimates of operators $\g\pa_t+\pa_x^4$ ($\g\in\dbR$). More precisely, we will establish two kinds of Carleman estimates with different weight functions, i.e.  singular weight function and regular weight function, respectively. Based on  Carleman estimate with singular weight function,  one can obtain the known controllability and observability results for the 1-d fourth order parabolic-type equation, while based on Carleman estimate with regular weight function, one can deduce not only the known result on conditional stability in the inverse problem of half-order parabolic equation, but also a new result on  H\"older type conditional stability in the inverse problem of half-order Schr\"odinger equation. We firmly believe that the Matlab program method developed in this paper can also be used to deal with other kind of Carleman estimate for the higher order partial differential operators.

Before explaining the key idea of Matlab program method, we first recall the key idea of Carleman estimate (see \cite{FLZ} for example).   Generally speaking, the main idea of Carleman estimate  is to multiply the principal part of partial differential operator by an appropriate weighted  multiplier. Then, the product is transformed into a sum of ``divergence" terms, ``energy" terms, and some low order  terms. Finally, by integrating this identity and choosing parameters large enough to absorb the undesired terms, which is the key idea of Carleman estimate.

For convenience, we give a simple example to introduce the key idea of Carleman estimate. Meanwhile, we give our  Matlab program method to see how to  compute it.

 { \bf Example:  } Let $\O=(0,1)$, for any fixed  $p\in
L^{\infty}( \O)$,  consider the
following second order operator:

$$
\cP(x,D)= \frac{d}{d x^2}+p(x),\q x\in \O.
$$

For any $x_0\in\dbR^n\setminus\cl{\O}$, set
 \bel{1pde2}
\Phi(x)=|x-x_0|^2.
 \ee

 \bp\label{Lemma1.2}  There exist two constants
$\l^*>0$ and $C>0$ such that
      \bel{1pde4}
\l^3\int_{\O}u^2e^{2\l\Phi}dx+\l\int_{\O}u_x^2e^{2\l\Phi}dx\le
C\int_{\O}|e^{\l\Phi}\cP(x,D)u|^2dx,\quad\forall\,\l
\ge\l^*, \ u\in H^2_0(\O).
   \ee
 \ep

{\it Proof. }
For $\l>0$, we set \bel{ell} \ell(x)=\l
\Phi(x),\q \th=\th(x)=e^{\ell(x)}, \q v=\th u, \ee
where $\Phi$ is given by (\ref{1pde2}). It is
easy to see that
  $$\ba{ll}\ds
\th
\cP(x,D)u&\ds=\th\big[(\th^{-1}v)_{xx}+p\th^{-1}v\big]=v_{xx}-2\ell_xv_x+(\ell_x^2-\ell_{xx})v+p v.  \ea
  $$
Put
$$
I_1=-2\ell_xv_x, \ I_2=v_{xx}+\ell_x^2 v-\ell_{xx}v, \ I_3=p v.
$$
 Then,
\bel{fst}\ba{ll}\ds
 \th \cP(x,D)u\cdot I_1\ds=I_1^2+I_1I_2+I_1I_3.
   \ea\ee
  Next, by (\ref{1pde2})  and (\ref{ell}), we have
   \begin{eqnarray*}
\begin{array}{rl}
 \ell_x=\l\mu, \q \mu=2(x-x_0), \q \ell_{xx}=2\l.
\end{array}
\end{eqnarray*}
A short calculation gives
  \bel{snd}\ba{ll}\ds
I_1I_2&\ds =-2\l\mu v_xv_{xx}-2\l^3\mu^3 vv_x+4\l^2\mu vv_x\\
    \ns&\ds=-(\l\mu v_x^2)_x+2\l v_x^2-(\l^3\mu^3v^2)_x+6\l^3\mu^2 v^2+(2\l^2\mu v^2)_x-4\l^2v^2.
  \ea\ee
Now, combining \eqref{fst} and \eqref{snd}, we have
  \bel{trd}\ba{ll}\ds
 \th \cP(x,D)u\cdot I_1\\
   \ns\ds =I_1^2+\big(-\l\mu v_x^2-\l^3\mu^3v^2+2\l^2\mu v^2)_x+2\l v_x^2+6\l^3\mu^2 v^2-4\l^2v^2-2\l\mu pvv_x.
  \ea\ee
By \eqref{trd}, we find that, for $c_0=4|p|^2_{L^\infty( \O)}\max\{(1-x_0)^2, x_0^2\}+4$,
$$
|\th \cP(x,D)u|^2 \ge
\big(-\l\mu v_x^2-\l^3\mu^3v^2+2\l^2\mu v^2)_x+2\l v_x^2+6\l^3\mu^2 v^2-c_0\l^2v^2-v_x^2.
$$
Finally,  integrating the above inequality in
$\O$, noting that $u\in H_0^2(\O)$ and $v=\th
u$, choosing $\l^*>0$ large enough,
then for all $\l\geq\l^*$, there exists a constant $C>0$, such that
inequality \eqref{1pde4} holds.
\endpf

Note that the key point in the proof of Proposition \ref{Lemma1.2} based on the pointwise identity  \eqref{trd}. To this aim, estimation of $I_1I_2$ is the key step.  In what follows, we  introduce the basic idea of using Matlab program method to calculate $I_1I_2$.

The key point of our method is to set a bijection to change all variable parameters into a vector (matrix).
For example, we define a bijection from $f(x)=a_0+a_1x+a_2x^2$ to a vector $\a=[a_0,a_1,a_2]$. Every element in vector $\a$ means the coefficient of $x^j$ ($j=0, 1, 2$).

Similarly, we define a bijection:
\bel{bi-01}\ba{ll}\ds
a_1 \l ^{a_2}[\mu^{a_3} v_{x^{(a_4)}}v_{x^{(a_5)}}]_{x^{(a_6)}}  \Leftrightarrow A=[a_1, a_2, a_3, a_4, a_5, a_6],
\ea\ee
where $\mu=2(x-x_0)$, $a_1$ is the constant coefficient, $a_2, a_3$ are the powers of $\l$, $\mu$, respectively, and $a_4, a_5, a_6$ are the derivative orders of $x$.  We observe that each term of $I_1I_2$ is of the form in \eqref{bi-01}, and $a_6$ is  always equal to 0.

Without loss of generality, we assume $a_4<a_5$.   Then, the computation of integration by parts can be replaced by number operation in matrix A.

If $a_4=a_5-1$, then,  by $\mu_x=2$, it is easy to see that
\bel{inte-001}\ba{ll}\ds
a_1 \l^{a_2}\mu^{a_3} v_{x^{(a_4)}}v_{x^{(a_5)}}=\frac{a_1}{2}\l^{a_2}[\mu^{a_3} v_{x^{(a_4)}}v_{x^{(a_5-1)}}]_{x}-a_1a_3\l^{a_2}\mu^{a_3-1} v_{x^{(a_4)}}v_{x^{(a_5-1)}}.
\ea\ee
The corresponding matrix operation is
\bel{oper-001}\ba{ll}\ds
[a_1, a_2, a_3, a_4, a_5, 0]\rightarrow \begin{bmatrix}
\frac{a_1}{2} & a_2 &a_3 & a_4 & a_5-1 & 1 \\
-a_1a_3 & a_2 & a_3-1 & a_4 & a_5-1 & 0
\end{bmatrix}.
\ea\ee
In this way, we can use matrix operations (\ref{oper-001})  to replace (\ref{inte-001}). First, write $I_1I_2$ in the form of a matrix
$$
I_1I_2=-2\l\mu v_xv_{xx}-2\l^3\mu^3 vv_x+4\l^2\mu vv_x \Leftrightarrow A=\begin{bmatrix}
-2& 1&1& 1 & 2& 0 \\
-2 &3& 3 & 0 & 1 & 0\\
4 &2& 1 & 0 & 1 & 0
\end{bmatrix}.
$$
By matrix operation in (\ref{oper-001}), the corresponding matrix becomes
\bel{oper-002}\ba{ll}\ds
 A=\begin{bmatrix}
-1 & 1 &1 & 1 & 1 & 1 \\
2 & 1 &0 & 1 & 1 & 0 \\
-1 & 3 & 3 & 0 & 0 & 1\\
-6 & 3 & 2& 0& 0 & 0\\
2 & 2 & 1 & 0 & 0 & 1\\
 -4& 2 & 0& 0& 0 & 0\\
\end{bmatrix}.
\ea\ee
By (\ref{bi-01}) we know that the translation result of  (\ref{oper-002}) is exactly the same as that of (\ref{snd}).

 This example is mainly to illustrate our ideas. In fact, we usually can't obtain the final result in just one loop. For second-order operators, we may not see the benefit of using matrix operations. However, for higher-order operators, there may be multiple matrix operations for each term, which can be seen in Section $2$. Therefore, as long as we define the matrix operations, the computer will automatically output the final result that meet our requirements.

Now, let us introduce  our main results. At the first beginning, we introduce some notations we used in this paper.
  Given $T>0$ and $L>0$, put
 $ \O=(0,L)$, $Q=(0,T)\t\O$.
 Throughout of this paper, we use $v_t$ or $\pa _t v$ to represent the derivative of $v$ in time, and the derivative of space also uses a similar notation.
 In what follows,  we will use $C$ to denote a generic
positive constant which may vary from line to line. The first Carleman estimate we obtained has the following form:

 \bt\label{thm-0} Let $\g\in\dbR$,  and put
  \bel{0a3-1}\ba{ll}\ds
 \psi(t,x)=(x-x_0)^2-\b(t-t_0)^2,\q \b>0,
 \ea\ee
 where $x_0>L$ such that $x_0-L>0$ is sufficiently small so that
$(x_0-L)^2<\frac{1}{3}x_0^2$ and $t_0 \in(0, T)$ be arbitrarily fixed. For any $\e>0$ satisfying $(x_0-L)^2<\e<x_0^2$, denote
 \bel{0a3x-1}\ba{ll}\ds
 $$
 Q_\e=\{(t,x)\in Q; \ \psi(t,x)>\e\}.
 $$
  \ea\ee
 Then for any $c\in (0,\e)$ and $v\in H^{(1,4)}(Q_c)$ with $\supp v\subset Q_\e$, there exists a $\l_0>0$, such that for any $\l>\l_0$, the following inequality holds:
  \bel{3a2-1}\ba{ll}\ds
  \int_{Q_{c}}e^{2\l\psi}\( \l^7|v|^2+\l^5| \pa _x v|^2+\l^3|\pa _{x}^2 v|^2+\l|\pa _x^3v|^2\)dxdt\\
  \ns\ds\le C\int_{Q_{c}}e^{2\l\psi} | \g\pa_t v+\pa_x^4 v|^2dxdt.
  \ea\ee
 \et

As application of Theorem \ref{thm-0}, we will discuss the conditional stability in inverse problem of a half-order fractional Schr\"odinger  equation (see Section \ref{section4} in detail).

 \br
In the case of $\g=-1$, inequality (\ref{3a2-1}) is given by \cite[Lemma 2]{YZ}. Based on this inequality, one can obtain the known result on  conditional stability in determining a zeroth-order coefficient in a half-order fractional diffusion equation.
 \er

 To obtain Carleman estimate with singular weight function for parabolic-type operator, we first recall give the choice of weight functions. Let $\o_0$ and $\o$ are two any given nonempty open subsets of $\O$ such that $\cl{\o_0} \subseteq \o$, by \cite{FI} we know that there exists a real-valued function $\ds\phi\in C^4(\cl{\O})$ such that
  \bel{0929-a}
  \phi> 0\q \mbox{ in }\O, \qq\phi=0\q\mbox{ on }\pa\O,\qq|\phi_x|>0\q \forall x\in \overline{\O\backslash \omega_0}.
  \ee

 For any $s>0$, put
  \bel{1a02}
  \rho(t ,x)=\frac{e^{s\phi(x)}-e^{2s|\phi|_{C(\cl{\O})}}}{t(T-t)}, \q \varphi(t, x)=\frac{e^{s\phi(x)}}{t(T-t)}.
 \ee

We have the following result.

  \bt\label{thm-01} Let $\rho, \varphi$ be given in (\ref{1a02}) and $\g\in\dbR$. Then there exists a positive constant $s_0$,  such that for any $s \ge s_0$, one can find two positive constants $\l_0=\l_0(\mu)$ and $C = C(s)$,  such that for all $v\in H^4(Q)$ satisfying
  $(\g\pa_t +\pa_x^4) v(t,x)=f$
  with $f \in L^2(Q)$ and $v=v_x=0$ on the boundary or $v=v_{xx}=0$ on the boundary, for all $\l\ge \l_0$, it holds that
  \bel{3a2-01}\ba{ll}\ds
 &\ds  \int_{Q}  \l s^2\varphi e^{2\l\rho} \( |\pa _x ^3v|^2+\l^2 s^2\varphi^2 |\pa_x^2 v|^2+\l^4 s^4\varphi ^4|\pa_x v|^2+\l^6 s^6\varphi^6 |v|^2\)dxdt\\
\ns&\ds\le C \int_{Q}e^{2\l\rho} f^2 dxdt+
C\int_{0}^T\int_{\omega}   \l^7 s^8\varphi^7e^{2\l\rho} |v|^2dxdt.
\ea\ee
 \et

 \br
Based on this Carleman estimate, we can obtain  the null controllability result of the fourth order parabolic-type equation (see \cite{ZZ} for example).
 \er

 The rest of this paper is organized as follows. In Section \ref{section2}, we will give the proof of Theorem \ref{thm-0}. In Section \ref{section3}, we will give the proof of Theorem \ref{thm-01}. In Section \ref{section4}, we will prove a conditional stability result in inverse source problem for a half-order time fractional Schr\"odinger equation.

 \section{Proof of  Theorem \ref{thm-0}}\label{section2}

In this Section, we will give the proof of Theorem \ref{thm-0}.

 {\bf Proof of Theorem \ref{thm-0}.} We divide the proof into several steps.

\ms

{\bf Step 1. }
To begin with, we  set
 \bel{1206-00b}
  \mu(x)=2(x-x_0),\q\cP\=\g\pa_t+\pa_x^4,
  \ee
and
 $$
\th=e^\ell, \q \ell=\l \psi, \q w=\th v,
$$
where $\psi$ is given by (\ref{0a3-1}).
Then it is easy to show that
 $$
 \ell_x=\l\mu(x)\=\l\mu,\q\ell_{xx}=2\l,\q\ell_{tx}=\ell_{xxx}=\ell_{xxxx}=0,
 $$
 and
 \bel{a1}\left\{\ba{ll}\ds
 \th v_t=w_t-\ell_tw,\q \th v_x=w_x-\l\mu w,\\
 \ns\ds \th v_{xx}= w_{xx}-2\l\mu w_x+(\l^2\mu^2-2\l)w.
 \ea\right.\ee
By (\ref{a1}), we have
 \bel{a2}\ba{ll}
&\ds\th v_{xxxx}=\[w_{xx}-2\l\mu w_x+(\l^2\mu^2-2\l)w\]_{xx}-2\l\mu\[w_{xx}-2\l\mu w_x+(\l^2\mu^2-2\l)w\]_x\\
\ns&\ds\q\q\q\q+(\l^2\mu^2-2\l)\[w_{xx}-2\l\mu w_x+(\l^2\mu^2-2\l)w\]\\
\ns&\ds\q\q\q=w_{xxxx}-4\l\mu w_{xxx}+(6\l^2\mu^2-12\l)w_{xx}+(-4\l^3\mu^3+24\l^2\mu)w_x\\
\ns&\ds\q\q\q\q+(\l^4\mu^4-12\l^3\mu^2+12\l^2)w.
 \ea\ee
Combining (\ref{a1}) and (\ref{a2}), we get
 \bel{a4}\ba{ll}\ds
\th(\g\pa_t v+\pa_x^4 v)=\th \cP v=I_1+I_2+I_3+I_4,
 \ea\ee
where
 \bel{a5}\left\{\ba{ll}\ds
  I_1\=\g w_t-4\l\mu
w_{xxx}-4\l^3\mu^3 w_x-12\l^3\mu^2 w,\\
\ns\ds I_2\=w_{xxxx}+6\l^2\mu^2w_{xx}+\l^4\mu^4w+12\l^2 w,\\
\ns\ds I_3\=-12\l w_{xx}+24\l^2\mu w_x,\\
\ns\ds I_4\=-\g\ell_tw.
 \ea\right.\ee

 In the case of $v$ is complex-valued function, by (\ref{a4}), we will compute
  \bel{za6}
\Re\int_{Q_{c}}\th\cP v\cl{I_1}dxdt=\Re\int_{Q_{c}}\(|I_1|^2+\cl{I_1}I_2+\cl{I_1}I_3+\cl{I_1}I_4\)dxdt.
 \ee
Noting that there is no $i=\sqrt{-1}$ in our principle operator $\cP$ and $I_j (j=1, 2, 3, 4)$,  for simplicity, we assume that $v$ is real-valued function, hence $w$ is real-valued function. In this situation, (\ref{za6}) can be rewritten as
 \bel{a6}
\int_{Q_{c}}\th\cP v I_1dxdt=\int_{Q_{c}}\(|I_1|^2+I_1I_2+I_1I_3+I_1I_4\)dxdt.
 \ee

{\bf Step 2. }  Since $I_1I_4$ is a low order term, the rest work is to  calculate $I_1I_2+I_1I_3$.
We note that each term of $I_1I_2+I_1I_3$ will be of the following form:
\bel{car-a1}\ba{ll}\ds
a_1 \l^{a_2}[\mu^{a_3}w_{t^{(a_4)}x^{(a_5)}}w_{x^{(a_6)}}]_{t^{(a_7)}x^{(a_8)}}  \Leftrightarrow A=[a_1,a_2,a_3,a_4,a_5, a_6, a_7, a_8],
\ea\ee
where $\mu$ is given by (\ref{1206-00b}), $a_1$ is the constant coefficient, $a_2, a_3$ are the powers of $\l$, $\mu$, respectively, $a_4,  a_7$ are the derivative orders of $t$,
and $a_5, a_6, a_8$ are the derivative orders of $x$.

Without loss of generality, we assume $a_5 \le a_6$. In fact, if $a_6 <a_5$, since $w$ is a real-valued function and $a_4$ must equal to $0$, we can just swap the positions of $a_5$ and $a_6$. Then, the computation of $I_1I_2+I_1I_3$ can be replaced by number operation in matrix $A$. First, we will explain how to convert the computation of integration by parts into operation in matrix $A$. Noting that  $a_4$ can only be equal to $1$ or $0$ in each term of $I_1I_2+I_1I_3$ and our purpose is to decompose $a_1\l^{a_2}\mu^{a_3} w_{t^{(a_4)}x^{(a_5)}}w_{x^{(a_6)}}$ into the sum of the divergence terms and the energy terms. To this aim, we divide our discussion into the following four cases.

If $a_7=1$ or $a_8=1$ or $a_4+a_5=a_6$, we skip this term. In other word, they are termination conditions.

If  $a_5=a_6$ and $a_4=1$,
then, by
\bel{inte-a0}\ba{ll}\ds
a_1 \l^{a_2}\mu^{a_3}w_{t^{(a_4)}x^{(a_5)}}w_{x^{(a_6)}}=\frac{a_1}{2}\l^{a_2}[\mu^{a_3} w_{t^{(a_4-1)}x^{(a_5)}}w_{x^{(a_6)}}]_{t},
\ea\ee
we can get the corresponding matrix operation as follows:
\bel{oper-a0}\ba{ll}\ds
\begin{bmatrix}
a_1 & a_2 & a_3 & a_4 & a_5 & a_6 & 0 & 0
\end{bmatrix}
\rightarrow \begin{bmatrix}
\frac{a_1}{2} & a_2 & a_3 & a_4-1 & a_5 & a_6 & 1 & 0
\end{bmatrix}.
\ea\ee

 If  $a_5=a_6-1$ and $a_4=0$,  noting $\mu_x=2$, it's easy to see that
\bel{inte-a1}\ba{ll}\ds
a_1 \l^{a_2}\mu^{a_3}w_{t^{(a_4)}x^{(a_5)}}w_{x^{(a_6)}}=\frac{a_1}{2}\l^{a_2}[\mu^{a_3} w_{t^{(a_4)}x^{(a_5)}}w_{x^{(a_6-1)}}]_{x}-a_1a_3\l^{a_2}\mu^{a_3-1} w_{t^{(a_4)}x^{(a_5)}}w_{x^{(a_6-1)}}.
\ea\ee
The corresponding matrix operation is
\bel{oper-a1}\ba{ll}\ds
[a_1, a_2, a_3, a_4, a_5, a_6, 0, 0]\rightarrow \begin{bmatrix}
\frac{a_1}{2} & a_2 & a_3 & a_4 & a_5 & a_6-1 & 0 & 1 \\
-a_1a_3 & a_2 & a_3-1 & a_4 & a_5 & a_6-1 & 0 & 0
\end{bmatrix}.
\ea\ee

 If $a_5< a_6-1$, by
 \bel{inte-a2}\ba{ll}\ds
a_1 \l^{a_2}\mu^{a_3}w_{t^{(a_4)}x^{(a_5)}}w_{x^{(a_6)}}=a_1[\l^{a_2}\mu^{a_3}w_{t^{(a_4)}x^{(a_5)}}w_{x^{(a_6-1)}}]_x-a_1\l^{a_2}\mu^{a_3}w_{t^{(a_4)}x^{(a_5+1)}}w_{x^{(a_6-1)}}\\
\ns\ds \q\q\q\q\q\q\q\q\q\q\q-2a_1a_3\l^{a_2}\mu^{a_3-1}w_{t^{(a_4)}x^{(a_5)}}w_{x^{(a_6-1)}},
\ea\ee
the corresponding matrix operation can be obtained as
\bel{oper-a2}\ba{ll}\ds
[a_1, a_2, a_3, a_4, a_5, a_6, 0, 0]\rightarrow \begin{bmatrix}
a_1 & a_2 & a_3 & a_4 & a_5 & a_6-1 & 0 & 1 \\
-a_1 & a_2 & a_3 & a_4 & a_5+1 & a_6-1 & 0 & 0 \\
-2a_1a_3 & a_2 & a_3-1 & a_4 & a_5 & a_6-1 & 0 & 0
\end{bmatrix}.
\ea\ee
In this way, we can use matrix operations (\ref{oper-a0}),  (\ref{oper-a1}) and (\ref{oper-a2}) to replace (\ref{inte-a0}), (\ref{inte-a1}) and (\ref{inte-a2}), respectively.

Next, we use one term in $I_1I_2+I_1I_3$ as an example to explain how this method works.
 First, just like (\ref{car-a1}), we set a matrix
\bel{210417-002}\ba{ll}\ds
-4\l^5\mu^5 ww_{xxx}\Leftrightarrow A= [-4, 5, 5, 0, 0, 3, 0, 0].
\ea\ee
This is equivalent to we first input the following  parameters:
 \begin{eqnarray*}\ba{ll}\ds
a_1=-4,\q a_2=5,\q a_3=5,\q a_4=0, \q a_5=0, \q a_6=3,\q a_7=0, \q a_8=0.
\ea\end{eqnarray*}
In the program,  we always default $a_7$ and $a_8$ to $0$ and we need not input it.

Because matrix A does not satisfy the termination conditions, we begin the 1st loop.
Since   $a_5< a_6-1$ in (\ref{210417-002}) , by (\ref{oper-a2}),  the corresponding matrix becomes
\bel{210417-04}\ba{ll}\ds
A=\begin{bmatrix}
-4 & 5 & 5 & 0 & 0 & 2 & 0 & 1 \\
4 & 5 & 5 & 0 & 1 & 2 & 0 & 0 \\
40 & 5 &4 & 0 & 0 & 2 &0 & 0
\end{bmatrix}.
\ea\ee

In the 2nd loop, we check each row of matrix A  in (\ref{210417-04}) whether it satisfies the termination condition.
The 1st row of matrix A meets the termination condition $a_8=1$, so we do nothing.
The 2nd row $[4, 5, 5, 0, 1, 2, 0, 0]$ of matrix A in (\ref{210417-04}) does not  meet the termination conditions,  so we apply integration by parts to this row. Notice that
  here  $a_5= a_6-1$ and $a_4=0$.  Then by (\ref{oper-a1}), we have
 \begin{eqnarray*}\ba{ll}\ds
[4, 5, 5, 0, 1, 2, 0, 0]\rightarrow \begin{bmatrix}
2 & 5 & 5 & 0& 1 & 1 & 0 & 1 \\
-20& 5 & 4 & 0 & 1 & 1 & 0 &0
\end{bmatrix}.
\ea\end{eqnarray*}
We noticed that the 3rd row $[40, 5, 4, 0, 0, 2, 0, 0]$ of matrix A in (\ref{210417-04}) also  does not meet the termination conditions.
By (\ref{oper-a2}), we can get that
\begin{eqnarray*}\ba{ll}\ds
[40, 5, 4, 0, 0, 2, 0, 0]\rightarrow \begin{bmatrix}
40 & 5 & 4 & 0 & 0 & 1 &0 & 1 \\
-40 & 5 & 4 & 0 & 1 & 1 &0 & 0 \\
-320& 5 & 3 & 0 & 0 & 1 &0 & 0
\end{bmatrix}.
\ea\end{eqnarray*}
Then the matrix A will become
 \bel{210417-08}\ba{ll}\ds
A=\begin{bmatrix}
-4 & 5 & 5 & 0 & 0 & 2 & 0 & 1 \\
2 & 5 & 5 & 0& 1 & 1 & 0 & 1  \\
40 & 5 & 4 & 0 & 0 & 1 &0 & 1 \\
-20& 5 & 4 & 0 & 1 & 1 & 0 &0\\
-40 & 5 & 4 & 0 & 1 & 1 &0 & 0 \\
-320& 5 & 3 & 0 & 0 & 1 &0 & 0
\end{bmatrix}.
\ea\ee

In loop 3rd,  we check all $6$ rows of $A$ in (\ref{210417-08}), and only the 6th row should be processed:
\begin{eqnarray*}\ba{ll}\ds
[-320, 5, 3, 0, 0, 1, 0, 0]\rightarrow \begin{bmatrix}
-160 & 5& 3 & 0 & 0 & 0& 0& 1 \\
960 & 5 & 2 & 0 & 0 & 0 &0 &0
\end{bmatrix}.
\ea\end{eqnarray*}
Then $A$ renew to
\bel{210417-10}\ba{ll}\ds
A=\begin{bmatrix}
-4 & 5 & 5 & 0 & 0 & 2 & 0 & 1 \\
2 & 5 & 5 & 0& 1 & 1 & 0 & 1  \\
40 & 5 & 4 & 0 & 0 & 1 &0 & 1 \\
-20& 5 & 4 & 0 & 1 & 1 & 0 &0\\
-40 & 5 & 4 & 0 & 1 & 1 &0 & 0 \\
-160 & 5& 3 & 0 & 0 & 0& 0& 1 \\
960 & 5 & 2 & 0 & 0 & 0 &0 &0
\end{bmatrix}.
\ea\ee

In loop 4th, every row of $A$ meets the termination conditions, so we output the matrix $A$.
On the other hand, by elementary calculations, it is easy to see that
\bel{210417-11}\ba{ll}\ds
-4\l^5\mu^5 ww_{xxx}&\ds=-4(\l^5\mu^5ww_{xx})_x+2(\l^5\mu^5w_x^2)_x+40(\l^5\mu^4ww_x)_x \\
\ns&\ds\q-60\l^5\mu^4w_x^2-160(\l^5\mu^3w^2)_x+960\l^5\mu^2w^2.
\ea\ee
Therefore, by (\ref{car-a1}), we can see that (\ref{210417-11}) is exactly the same as (\ref{210417-10}).

According to the above analysis, now we give the rule for calculating each term in $I_1I_2+I_1I_3$. It should be pointed out that in our code, our input data's format is $a_1 \lambda^{a_2} \mu^{a_3} w_{t^{(a_4)}x^{(a_5)}}$ from (\ref{a5}). Meanwhile, noting that $\g w_t I_2+\g w_tI_3$ is a low-order term, so we set $\g=1$ in the program, which will not affect the final result.

{\bf Algorithm 2.1. }

(1). Input matrix  $A$ with 8 columns. Loop (2) until all rows are skipped, then output $A$.

(2). Check every row $a\=[a_1,a_2,a_3,a_4,a_5,a_6,a_7,a_8]$ in $A$. If $a_1=0$ or $a_7=1$ or $a_8=1$ or $a_4+a_5=a_6$, then skip this row, else do (3).

(3). If $a_5=a_6$, then $a\rightarrow[a_1/2 ,a_2,a_3, a_4-1 ,a_5, a_6, 1, 0]$, else do (4).

(4). If $a_5> a_6$, then do (4.1), else do (4.2).

(4.1). If $a_4=0$ and $a_5-1=a_6$,  then
\begin{eqnarray*}\ds
a\rightarrow\begin{bmatrix}
a_1/2 & a_2 & a_3 & a_4 & a_5-1 & a_6 & 0 & 1 \\
-a_1a_3 & a_2 & a_3-1 & a_4 & a_5-1 & a_6 & 0 & 0
\end{bmatrix};
 \end{eqnarray*}
 else
 \begin{eqnarray*}\ds
a\rightarrow\begin{bmatrix}
a_1 & a_2 & a_3 & a_4 & a_5-1 & a_6 & 0 & 1 \\
-a_1 & a_2 & a_3 & a_4 & a_5-1 & a_6+1 & 0 & 0 \\
-2a_1a_3 & a_2 & a_3-1 & a_4 & a_5-1 & a_6 & 0 & 0
\end{bmatrix}.
 \end{eqnarray*}

(4.2). If $a_4=0$ and $a_6-1=a_5$,  then
\begin{eqnarray*}\ds
a\rightarrow\begin{bmatrix}
a_1/2 & a_2 & a_3 & a_4 & a_5 & a_6-1 & 0 & 1 \\
-a_1a_3 & a_2 & a_3-1 & a_4 & a_5 & a_6 -1& 0 & 0
\end{bmatrix};
 \end{eqnarray*}
 else
\begin{eqnarray*}\ds
a\rightarrow\begin{bmatrix}
a_1& a_2 & a_3 & a_4 & a_5 & a_6-1 & 0 & 1 \\
-a_1 & a_2 & a_3 & a_4 & a_5+1 & a_6-1 & 0 & 0 \\
-2a_1a_3 & a_2 & a_3-1 & a_4 & a_5 & a_6 -1& 0 & 0
\end{bmatrix}.
 \end{eqnarray*}

 {\bf Step 3. }
Noting that $c\in(0,\e)$ and $\supp v\subset Q_{\e}$, then $v$ equals to zero on both the space boundary and the time boundary of $Q_c$. Therefore, we delete the boundary term $a_7=1$ or $a_8=1$ in  matrix A, and extract the highest order energy terms at the same time. If the highest order energy terms is positive, we delete the lower order terms, otherwise we keep the lower order terms.
 Use the above method to estimate $I_1I_2+I_1I_3$, we have the following output result:
 \begin{figure}[H]
 \centering
 \includegraphics[width=0.80\textwidth]{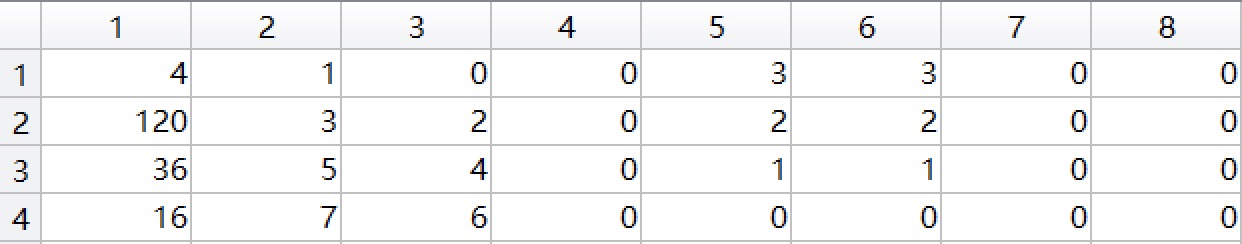}
 \label{01}
 \caption{the highest order energy terms of $I_1I_2+I_1I_3$}
 \end{figure}

From Figure 1, we can get, for sufficiently large $\l$,
 \begin{eqnarray}\label{i1i2}
\begin{array}{rl}
&\ds \int_{Q_{c}}\(I_1I_2+I_1I_3\)dxdt\\
\ns&\ds\geq \int_{Q_{c}} \(4\l w_{xxx}^2+120\l^3\mu^2w_{xx}^2+36\l^5\mu^4w_x^2+16\l^7\mu^6w^2\)dxdt.
\end{array}\end{eqnarray}
Combining(\ref{a5}),  (\ref{a6}) and (\ref{i1i2}), we conclude that there exists a positive constant $c_0\in (0, 1)$,  for sufficiently large $\l$, it holds that
\begin{eqnarray*}
\begin{array}{rl}
&\ds \int_{Q_{c}}\th\cP v I_1dxdt=\int_{Q_{\e_0}}\(|I_1|^2+I_1I_2+I_1I_3+I_1I_4\)dxdt\\
\ns&\ds\geq \int_{Q_{c}} \(c_0|I_1|^2+4\l w_{xxx}^2+120\l^3\mu^2w_{xx}^2+36\l^5\mu^4w_x^2+16\l^7\mu^6w^2\)dxdt.
\end{array}\end{eqnarray*}
 Notice that $w=\th v$, we can easily obtain (\ref{3a2-1}).\endpf

 \section{Proof of  Theorem \ref{thm-01}}\label{section3}

 In this Section, we will give the proof of Theorem \ref{thm-01}.

  {\bf Proof of Theorem \ref{thm-01}.} We divide the proof into several steps.

{\bf Step 1. }
To begin with, we  set
\bel{1206-0b0}
\th=e^\ell, \q \ell=\l\rho,  \q u=\th v,
\ee
where $\rho$ is given by (\ref{1a02}).
By  (\ref{1206-0b0}), it is easy to check
that
 \bel{b1}\left\{\ba{ll}\ds
 \th v_t=u_t-\ell_tu,\q \th v_x=u_x-\ell_x u,\\
 \ns\ds \th v_{xx}= u_{xx}-2\ell_x u_x+(\ell_x^2-\ell_{xx})u.
 \ea\right.\ee
Similar to (\ref{a2}), we can get
 \bel{b2}\ba{ll}
&\th v_{xxxx}\ds
=u_{xxxx}-4\ell_x u_{xxx}+6(\ell_x^2-\ell_{xx})u_{xx}+(-4\ell_x^3+12\ell_x\ell_{xx}-4\ell_{xxx})u_x\\
\ns&\ds\q\q\qq+(\ell_x^4-6\ell_x^2 \ell_{xx}+3\ell_{xx}^2+4\ell_x\ell_{xxx}-\ell_{xxxx})u.
 \ea\ee
Therefore, combining (\ref{1206-0b0})--(\ref{b2}),  we get
 \bel{b4}\ba{ll}\ds
\th(\g\pa_t v+\pa_x^4 v)=\th \cP v=J_1+J_2+J_3+J_4,
 \ea\ee
where
 \bel{b5}\left\{\ba{ll}\ds
  J_1\=\g u_t-4\l\rho_x
u_{xxx}-4\l^3\rho_x^3 u_x-6\l^3\rho_x^2\rho_{xx} u,\\
\ns\ds J_2\=u_{xxxx}+6\l^2\rho_x^2u_{xx}+\l^4\rho_x^4u+(3\l^2 \rho_{xx}^2+4\l^2 \rho_{x}\rho_{xxx})u,\\
\ns\ds J_3\=-6\l \rho_{xx} u_{xx}+(12\l^2\rho_x \rho_{xx}-4\l \rho_{xxx}) u_x,\\
\ns\ds J_4\=(-\l \rho_{xxxx}-\g\ell_{t})u.
 \ea\right.\ee
For simplicity, we assume that $v$ is real-valued function, hence, by (\ref{b4}), we have
  \bel{zb6}
\th\cP vJ_1=|J_1|^2+J_1J_2+J_1J_3+J_1J_4.
 \ee

{\bf Step 2. }
Since $J_1J_4$ is a low order term, the rest work is to estimate $J_1J_2+J_1J_3$.
By the definitions of $\rho$ in (\ref{1a02}), we have that
  \bel{b3}\left\{\ba{ll}\ds
  \rho_x=s\phi_x\varphi,\q \rho_{xx}=s^2\phi_x^2\varphi+ s\phi_{xx}\varphi,\\
 \ns\ds  \rho_{xxx}=s^3\phi_x^3\varphi+3 s^2 \phi_x\phi_{xx}\varphi+ s \phi_{xxx}\varphi.
 \ea\right.\ee
Then we can see that each term of $J_1J_2+J_1J_3$ can be written as follows:
\bel{car-b1}\ba{ll}\ds
&a_1 \l^{a_2} s^{a_3} [\phi_x^{a_4} \phi_{xx}^{a_5}\phi_{xxx}^{a_6} \phi_{xxxx}^{a_7}\varphi^{a_8}\varphi_t^{a_9}  u_{t^{(a_{10})}x^{(a_{11})}}u_{x^{(a_{12})}}]_{t^{(a_{13})}x^{(a_{14})}}\\
 \ns&\ds \Leftrightarrow A=[a_1, a_2, a_3, a_4, a_5, a_6, a_7 , a_8, a_9, a_{10}, a_{11}, a_{12}, a_{13}, a_{14}].
\ea\ee
where $a_1$ is the constant coefficient, $a_2, a_3, a_4, a_5, a_6, a_7 , a_8, a_9$ are the powers of $\l$, $s$, $\phi_x$, $\phi_{xx}$, $\phi_{xxx}$, $\phi_{xxxx}$, $\varphi$, $\varphi_t$, respectively, $a_{10}$,  $a_{13}$ are the derivative orders of $t$,
and $a_{11}, a_{12}, a_{14}$ are the derivative orders of $x$.

The method for calculating $J_1J_2 + J_1J_3$ is similar to the idea in the previous subsection, except that  it requires a larger matrix to store the corresponding positions of $s$, $\phi_x$, $\phi_{xx}$, $\phi_{xxx}$, $\phi_{xxxx}$, $\varphi$ and $\varphi_t$. Next, we use one term in $J_1J_2+J_1J_3$ as an example to explain how our method works.

First, just like (\ref{car-b1}), we set a matrix
\bel{210417-b2}\ba{ll}\ds
-36\l^5 s^5\phi_x^5\varphi^4uu_{xx}
\Leftrightarrow
\addtocounter{MaxMatrixCols}{10}
A= \begin{bmatrix}
-36 & 5 & 5 & 5 & 0 & 0 & 0 & 4 & 0 & 0 & 0 & 2 & 0 & 0 \\
\end{bmatrix}.
\ea\ee

Similarly, we define termination conditions as follows: $a_{13}=1$ or $a_{14}=1$ or $a_{10}+a_{11}=a_{12}$.
Note that matrix A  in (\ref{210417-b2}) does not satisfy the termination condition, we begin the 1st loop using integration by parts.
Decompose matrix A for the 1st loop, and the corresponding matrix becomes:
\bel{210417-b4}\ba{ll}\ds
\addtocounter{MaxMatrixCols}{10}
A=\begin{bmatrix}
-36 & 5 & 5 & 5 & 0 & 0 & 0 & 4 & 0 & 0 & 0 & 1 & 0 & 1 \\
36 & 5 & 5 & 5 & 0 & 0 & 0 & 4 & 0 & 0 & 1 & 1 & 0 & 0 \\
180 & 5 & 5 & 4 & 1 & 0 & 0 & 4 & 0 & 0 & 0 & 1 & 0 & 0 \\
144 & 5 & 6 & 6 & 0 & 0 & 0 & 4 & 0 & 0 & 0 & 1 & 0 & 0
\end{bmatrix}.
\ea\ee

In the 2nd loop, we check each row of matrix A.
The 1st and  2nd row of matrix A meets the termination conditions, so we do nothing.
The 3rd row of matrix A in (\ref{210417-b4}) does not meet the termination conditions. Doing integration by parts, we get
  \begin{eqnarray*}\ba{ll}\ds
  \addtocounter{MaxMatrixCols}{10}
  \begin{bmatrix}
180 & 5 & 5 & 4 & 1 & 0 & 0 & 4 & 0 & 0 & 0 & 1 & 0 & 0
\end{bmatrix}\\
\ns\ds\rightarrow \begin{bmatrix}
90 & 5 & 5 & 4 & 1 & 0 & 0 & 4 & 0 & 0 & 0 & 0 & 0 & 1 \\
-360 & 5 & 5 & 3 & 2 & 0 & 0 & 4 & 0 & 0 & 0 & 0 & 0 & 0 \\
-90 & 5 & 5 & 4 & 0 & 1 & 0 & 4 & 0 & 0 & 0 & 0 & 0 & 0 \\
-360 & 5 & 6 & 5 & 1 & 0 & 0 & 4 & 0 & 0 & 0 & 0 & 0 & 0
\end{bmatrix}.
\ea\end{eqnarray*}
The 4th row of matrix A in (\ref{210417-b4}) also does not meet the termination conditions.
By
 \begin{eqnarray*}\ba{ll}\ds
  \addtocounter{MaxMatrixCols}{10}
  \begin{bmatrix}
144 & 5 & 6 & 6 & 0 & 0 & 0 & 4 & 0 & 0 & 0 & 1 & 0 & 0
\end{bmatrix}\\
\ns\ds\rightarrow \begin{bmatrix}
72 & 5 & 6 & 6 & 0 & 0 & 0 & 4 & 0 & 0 & 0 & 0 & 0 & 1 \\
-432 & 5 & 6 & 5 & 1 & 0 & 0 & 4 & 0 & 0 & 0 & 0 & 0 & 0 \\
-288 & 5 & 7 & 7 & 0 & 0 & 0 & 4 & 0 & 0 & 0 & 0 & 0 & 0
\end{bmatrix},
\ea\end{eqnarray*}
 then $A$ renew to
\bel{210417-b6}\ba{ll}\ds
\addtocounter{MaxMatrixCols}{10}
A=\begin{bmatrix}
-36 & 5 & 5 & 5 & 0 & 0 & 0 & 4 & 0 & 0 & 0 & 1 & 0 & 1 \\
36 & 5 & 5 & 5 & 0 & 0 & 0 & 4 & 0 & 0 & 1 & 1 & 0 & 0 \\
90 & 5 & 5 & 4 & 1 & 0 & 0 & 4 & 0 & 0 & 0 & 0 & 0 & 1 \\
72 & 5 & 6 & 6 & 0 & 0 & 0 & 4 & 0 & 0 & 0 & 0 & 0 & 1 \\
-360 & 5 & 5 & 3 & 2 & 0 & 0 & 4 & 0 & 0 & 0 & 0 & 0 & 0 \\
-90 & 5 & 5 & 4 & 0 & 1 & 0 & 4 & 0 & 0 & 0 & 0 & 0 & 0 \\
-360 & 5 & 6 & 5 & 1 & 0 & 0 & 4 & 0 & 0 & 0 & 0 & 0 & 0 \\
-432 & 5 & 6 & 5 & 1 & 0 & 0 & 4 & 0 & 0 & 0 & 0 & 0 & 0 \\
-288 & 5 & 7 & 7 & 0 & 0 & 0 & 4 & 0 & 0 & 0 & 0 & 0 & 0
\end{bmatrix}.
\ea\ee

In loop 3rd, since every row of $A$ meets the termination conditions, we output the matrix $A$.
The translation from (\ref{210417-b6}) means
\begin{eqnarray*}\ba{ll}\ds
&-36\l^5 s^5\phi_x^5\varphi^4uu_{xx}\\
 \ns&\ds=-36\l^5 s^5(\phi_x^5\varphi^4uu_x)_x+36\l^5 s^5\phi_x^5\varphi^4u_x^2+90\l^5 s^5(\phi_x^4\phi_{xx}\varphi^4u^2)_x+72\l^5 s^6(\phi_x^6\varphi^4u^2)_x\\
\ns&\ds
\q -360\l^5 s^5\phi_x^3\phi_{xx}^2\varphi^4u^2-90\l^5 s^5\phi_x^4\phi_{xxx}\varphi^4u^2-792\l^5 s^6\phi_x^5\phi_{xx}\varphi^4u^2-288\l^5 s^7\phi_x^7\varphi^4u^2.
\ea\end{eqnarray*}
According to the above analysis, we first write $J_1J_2+J_1J_3$ in the form of a matrix, and then obtain a matrix similar to (\ref{210417-b6}) through matrix operation. Finally, we extract the boundary terms, cross terms and energy terms,  respectively.

 For the energy terms, we keep the highest and subhigher order terms  of $\l$ and $s$, while deleting other lower order terms which can be absorbed by the  subhigher order terms. The results of the highest order energy terms and cross terms are shown in Figure 2, the results of the subhigher order terms  are shown in Figure 3.  In the following, for any $k \in \mathbb{N}$, we denote by  $\mathcal{O}(s^k)$  a function of order $s^k$ for  large $s$.
 \begin{figure}[H]
 \centering
 \includegraphics[width=16cm,height=3cm]{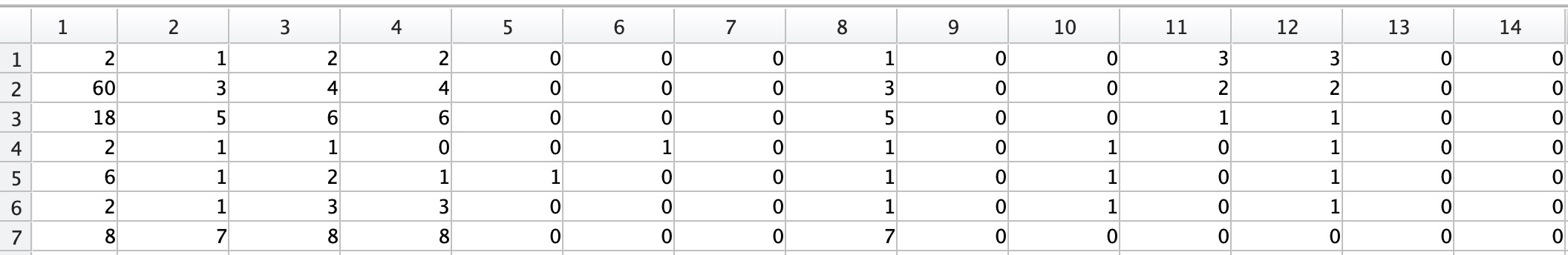}
 \label{02}
 \caption{the highest order energy terms and cross terms of $J_1J_2+J_1J_3$}
 \end{figure}

 \begin{figure}[H]
 \centering
 \includegraphics[width=16cm,height=3cm]{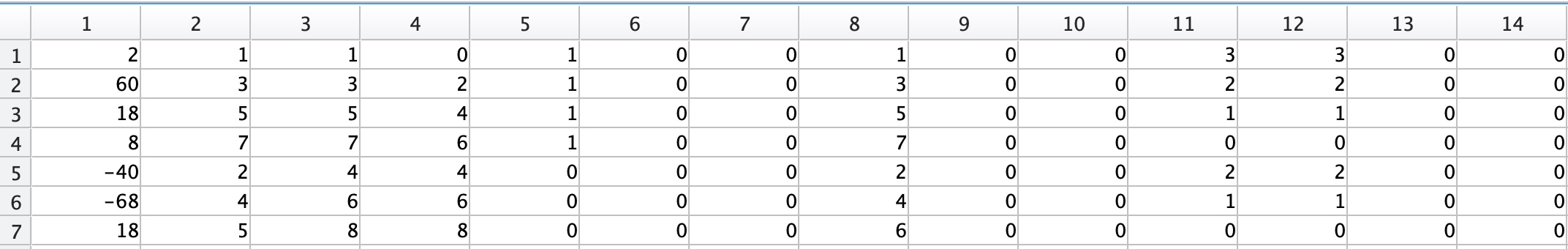}
 \label{03}
 \caption{subhigher order terms  of $J_1J_2+J_1J_3$}
 \end{figure}
For the time boundary, $\th$  is singular with respect to time variable,  we directly remove the rows satisfying $a_{13}=1$ in matrix A.
 For the space boundary, i.e. the row $a_{14}=1$ in  matrix $A$, if $v=v_x=0$ on the boundary, we delete the boundary terms with $u$ or $u_x$ or $u_t$ or $u_{tx}$ on the boundary terms, and the results of the boundary terms are shown in Figure 4.
  \begin{figure}[H]
  \centering
 \includegraphics[width=16cm,height=2cm]{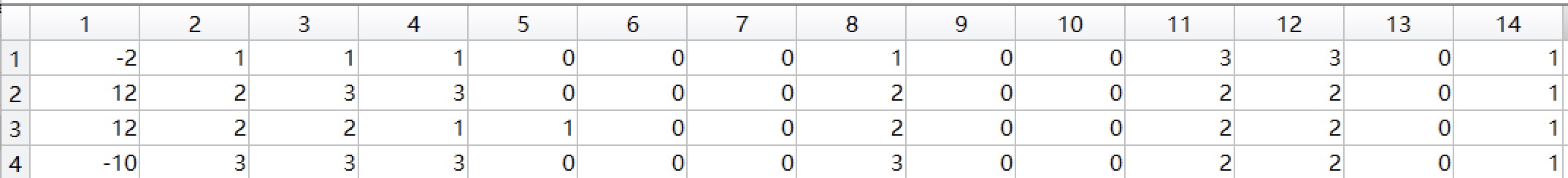}
 \label{04}
 \caption{the boundary terms of $J_1J_2+J_1J_3$}
 \end{figure}

 \br
The reason why the algorithm is not proposed here is that except for the different matrix operations, the calculation rules here are similar to Algorithm 2.1.
 \er

\ms
 {\bf Step 3. }
 From Figure 1, 2, 3, we can get
\begin{eqnarray}\label{21092901}
\begin{array}{rl}
&\ds \int_{Q}\(J_1J_2+J_1J_3\)dxdt\\
\ns&\ds=  \int_{Q} \(2\l s^2\varphi \phi_{x}^2 u_{xxx}^2+60\l^3s^4\varphi^3 \phi_{x}^4u_{xx}^2+18\l^5s^6\varphi ^5\phi_{x}^6u_x^2+8\l^7s^8\varphi^7 \phi_{x}^8u^2+Bu_tu_x\)dxdt\\
\ns&\ds\q + \int_{Q}V_xdxdt+\int_{Q} \Big\{\mathcal{O}(s)\l\varphi u_{xxx}^2+\big[\mathcal{O}(\l^3) s^3\varphi^3+\mathcal{O}(s^4) \l^2\varphi^2\big]u_{xx}^2\Big\}dxdt\\
\ns&\ds\q+\int_{Q}\Big\{\big[\mathcal{O}(\l^5) s^5\varphi^5+\mathcal{O}(s^4) \l^4\varphi^4\big]u_{x}^2+\big[\mathcal{O}(\l^7) s^7\varphi^7+\mathcal{O}(s^8) \l^5\varphi^6\big]u^2\Big\}dxdt,
\end{array}\end{eqnarray}
where
\ms
\begin{eqnarray}\label{21092902}
\begin{array}{rl}
&\ds V= -2\l s\varphi \phi_x u_{xxx}^2
-10\l^3 s^3\varphi^3 \phi_x^3 u_{xx}^2+\mathcal{O}(s^3)\l^2\varphi^2u_{xx}^2,
\end{array}\end{eqnarray}
and
\begin{eqnarray}\label{cross}
\begin{array}{rl}
&B= 2\l s \varphi  \phi_{xxx}+6\l s^2\varphi\phi_x\phi_{xx}+2\l s^3\varphi\phi_x^3.
\end{array}\end{eqnarray}

In the case of $v=v_{xx}=0$ on the boundary, we can similarly deal with it. Extract the terms $a_{14}=1$ from the matrix $A$ and delete the terms with $u$ or $u_t$ or $u_{xx}$ or $u_{txx}$,  leaving only the highest order terms, the  subhigher order terms and the cross terms at the end. Other low-order terms can be absorbed by the  subhigher order terms. The results of the boundary terms in this case are shown in Figure 5.
 \begin{figure}[H]
  \centering
 \includegraphics[width=16cm,height=2cm]{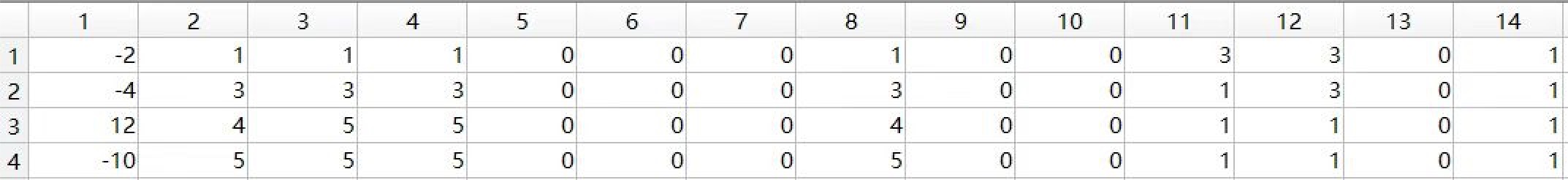}
 \label{05}
 \caption{the boundary terms }
 \end{figure}
 Figure 5 shows that the boundary term in this case is
\begin{eqnarray}\label{21092903}
\begin{array}{rl}
&\ds V= -2\l s\varphi \phi_x u_{xxx}^2
-4\l^3 s^3\varphi^3 \phi_x^3 u_xu_{xxx}-10\l^5 s^5\varphi^5\phi_x^5u_x^2+\mathcal{O}(s^5)\l^4\varphi^4u_{x}^2.
\end{array}\end{eqnarray}

{\bf Step 4. }
By the definition of $J_1$ in (\ref{b5}) and $B$ in (\ref{cross}),  noting that (\ref{b3}), we have, for some $\e>0$,
\begin{eqnarray}\label{c2}
\begin{array}{rl}
&\ds B u_t u_x=Bu_x\(J_1+4\l s \varphi \phi_xu_{xxx}+4\l ^3s^3 \varphi^3 \phi_x^3u_{x}+6\l^3 s^4\varphi^3\phi_x^4u+6\l^3 s^3\varphi^3\phi_x^2\phi_{xx}u\)\\
\ns&\ds\ge Bu_x\(4\l s \varphi \phi_xu_{xxx}+4\l ^3s^3 \varphi^3 \phi_x^3u_{x}+6\l^3 s^4\varphi^3\phi_x^4u+6\l^3 s^3\varphi^3\phi_x^2\phi_{xx}u\)\\
\ns&\ds\q -\e J_1^2-C\l^2s^6\varphi^2u_x^2.
\end{array}\end{eqnarray}
Notice that the right-hand side of inequality  (\ref{c2}) can be rewritten as
matrix with the same form as (\ref{car-b1}). Therefore, we can apply the method again. For  energy terms and boundary terms, only the highest order terms are left. The result is shown in Figure 6.
 \begin{figure}[H]
  \centering
 \includegraphics[width=16cm,height=2.8cm]{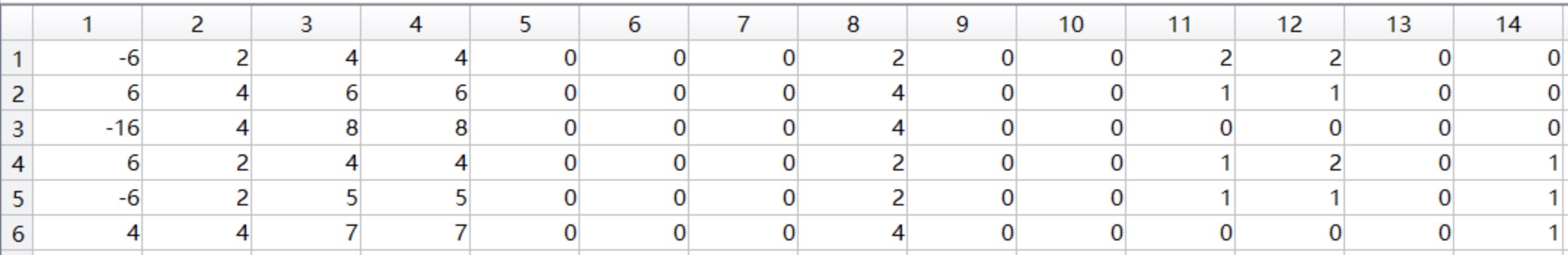}
 \label{06}
 \caption{the energy terms and boundary terms of $Bu_tu_x$ }
 \end{figure}
 Then, we have
\begin{eqnarray}\label{21092904}
\begin{array}{rl}
&\ds   \int_{Q} B u_t u_x dxdt\\
\ns&\ds\ge\int_{Q}  \(-\e J_1^2-C\l^2s^4\varphi^2u_{xx}^2-C\l^4s^6\varphi^4 u_{x}^2-C\l^4s^8\varphi^4u^2\)dxdt+ \int_{Q} \widetilde{V}_xdxdt.
\end{array}\end{eqnarray}
where
$$
 \widetilde{V}=\mathcal{O}(s^4)\l^2\varphi^2u_xu_{xx}+\mathcal{O}(s^4)\l^7\varphi^4u^2+\mathcal{O}(s^5)\l^2\varphi^2u_{x}^2
$$
Therefore, combining (\ref{21092901}) and (\ref{21092904}), we have
\begin{eqnarray*}
\begin{array}{rl}
&\ds \int_{Q} \(2\l s^2\varphi \phi_{x}^2 u_{xxx}^2+60\l^3s^4\varphi^3 \phi_{x}^4u_{xx}^2+18\l^5s^6\varphi ^5\phi_{x}^6u_x^2+8\l^7s^8\varphi^7 \phi_{x}^8u^2\)dxdt+ \int_{Q} V_xdxdt\\
\ns&\ds\le \int_{Q}\(J_1J_2+J_1J_3+\e J_1^2\)dxdt+C\int_{Q} \Big[s\l\varphi u_{xxx}^2+(\l^3 s^3\varphi^3+ \l^2s^4 \varphi^2)u_{xx}^2\Big]dxdt\\
\ns&\ds\q+C\int_{Q}\Big[(\l^5 s^5\varphi^5+ \l^4s^6 \varphi^4)u_{x}^2+(\l^7 s^7\varphi^7+ \l^5s^8\varphi^6)u^2\Big]dxdt.
\end{array}\end{eqnarray*}
where,  for $v=v_{x}=0$ and $v=v_{xx}=0$ on the boundary,  $V$ is given by (\ref{21092902}) and (\ref{21092903}),  respectively.
By (\ref{0929-a}), we know that $\phi_x(L)<0$ and $\phi_x(0)>0$,  then for sufficiently large $s$,  we have
\begin{eqnarray*}
\begin{array}{rl}
&\ds V(L,t)\ge 0, \  V(0,t)\le 0.
\end{array}\end{eqnarray*}
Hence,
\begin{eqnarray}\label{c4}
\begin{array}{rl}
&\ds \int_{Q} \(2\l s^2\varphi \phi_{x}^2 u_{xxx}^2+60\l^3s^4\varphi^3 \phi_{x}^4u_{xx}^2+18\l^5s^6\varphi ^5\phi_{x}^6u_x^2+8\l^7s^8\varphi^7 \phi_{x}^8u^2\)dxdt\\
\ns&\ds\le \int_{Q}\(J_1J_2+J_1J_3+\e J_1^2\)dxdt+C\int_{Q} \Big[s\l\varphi u_{xxx}^2+(\l^3 s^3\varphi^3+\l^2s^4 \varphi^2)u_{xx}^2\Big]dxdt\\
\ns&\ds\q+C\int_{Q}\Big[(\l^5 s^5\varphi^5+\l^4s^6 \varphi^4)u_{x}^2+(\l^7 s^7\varphi^7+ \l^5 s^8 \varphi^6)u^2\Big]dxdt.
\end{array}\end{eqnarray}
Then, it follows from (\ref{zb6}) and (\ref{c4}) that, for sufficiently small $\e>0$,
\begin{eqnarray}\label{c5}
\begin{array}{rl}
&\ds  \int_{Q}  \(\l s^2\varphi \phi_{x}^2 u_{xxx}^2+\l^3s^4\varphi^3 \phi_{x}^4u_{xx}^2+\l^5s^6\varphi ^5\phi_{x}^6u_x^2+\l^7s^8\varphi^7 \phi_{x}^8u^2\)dxdt\\
\ns&\ds\le C\int_{Q}\th^2|\cP v|^2 dxdt+C\int_{Q} \Big[s\l\varphi u_{xxx}^2+(\l^3 s^3\varphi^3+\l^2s^4 \varphi^2)u_{xx}^2\Big]dxdt\\
\ns&\ds\q+C\int_{Q}\Big[(\l^5 s^5\varphi^5+\l^4s^6 \varphi^4)u_{x}^2+(\l^7 s^7\varphi^7+\l^5s^8 \varphi^6)u^2\Big]dxdt.
\end{array}\end{eqnarray}
Further, noting that $|\phi_x| > 0$ in $\O\backslash \omega_0$,  it follows that
\begin{eqnarray}\label{c6}
\begin{array}{rl}
&\ds   \int_{Q}  \(\l s^2\varphi \phi_{x}^2 u_{xxx}^2+\l^3s^4\varphi^3 \phi_{x}^4u_{xx}^2+\l^5s^6\varphi ^5\phi_{x}^6u_x^2+\l^7s^8\varphi^7 \phi_{x}^8u^2\)dxdt\\
\ns&\ds\ge C\int_{0}^T\int_{\O\backslash\omega_0}  \(\l s^2\varphi  u_{xxx}^2+\l^3s^4\varphi^3 u_{xx}^2+\l^5s^6\varphi ^5u_x^2+\l^7s^8\varphi^7 u^2\)dxdt\\
\ns&\ds\q- C \int_{0}^T\int_{\omega_0}  \(\l s^2\varphi  u_{xxx}^2+\l^3s^4\varphi^3 u_{xx}^2+\l^5s^6\varphi ^5u_x^2+\l^7s^8\varphi^7 u^2\)dxdt.
\end{array}\end{eqnarray}
Consequently, it follows from (\ref{c6}) that
\begin{eqnarray}\label{c7}
\begin{array}{rl}
&\ds   \int_{Q}  \(\l s^2\varphi  u_{xxx}^2+\l^3s^4\varphi^3 u_{xx}^2+\l^5s^6\varphi ^5u_x^2+\l^7s^8\varphi^7u^2\)dxdt\\
\ns&\ds\le  C\int_{Q}  \(\l s^2\varphi \phi_{x}^2 u_{xxx}^2+\l^3s^4\varphi^3 \phi_{x}^4u_{xx}^2+\l^5s^6\varphi ^5\phi_{x}^6u_x^2+\l^7s^8\varphi^7 \phi_{x}^8u^2\)dxdt\\
\ns&\ds\q+C \int_{0}^T\int_{\omega_0}  \(\l s^2\varphi  u_{xxx}^2+\l^3s^4\varphi^3 u_{xx}^2+\l^5s^6\varphi ^5u_x^2+\l^7s^8\varphi^7 u^2\)dxdt.
\end{array}\end{eqnarray}
Combining (\ref{c5}) and (\ref{c7}), and noting $u=\th v$,  we conclude that there exists a $s_0>0$ such that for all $s\ge s_0$, one can find a positive constant $\l=\l_0(s)$ so that for any $\l\ge\l_0$, it holds that
\begin{eqnarray}\label{c8}
\begin{array}{rl}
&\ds  \int_{Q}  \th^2\l s^2\varphi \( v_{xxx}^2+\l^2s^2\varphi^2 v_{xx}^2+\l^4s^4\varphi ^4v_x^2+\l^6s^6\varphi^6 v^2\)dxdt\\
\ns&\ds\le C\int_{0}^T\int_{\omega_0}  \th^2\l s^2\varphi \( v_{xxx}^2+\l^2s^2\varphi^2 v_{xx}^2+\l^4s^4\varphi ^4v_x^2+\l^6s^6\varphi^6 v^2\)dxdt\\
\ns&\ds\q+ C \int_{Q}\th^2|\cP v|^2 dxdt.
\end{array}
\end{eqnarray}

{\bf Step 5. }
We shall eliminate the terms $\ds\int_{0}^T\int_{\omega_0}  \th^2\l s^2\varphi \( v_{xxx}^2+\l^2s^2\varphi^2 v_{xx}^2+\l^4s^4\varphi ^4v_x^2\)dxdt$.

Take $\zeta\in C_0^{\infty}(\omega; [0,1])$ such that $\zeta=1$ in $\omega_0$. Multiplying $\g v_t+v_{xxxx}$ by $\zeta^8\th^2\varphi v_{xx}$ and integrating it over $Q$, we get
\begin{eqnarray*}
\begin{array}{rl}
&\ds  \int_{Q} \[ (\g v_t+v_{xxxx})\zeta^8\th^2\varphi v_{xx}\]dxdt\\
\ns&\ds =\int_{Q}  \[-\g(\zeta^8\th^2\varphi)_x v_tv_x+\frac{1}{2}\g\zeta^8(\th^2\varphi)_tv_x^2- (\zeta^8\th^2\varphi)_{x}v_{xx}v_{xxx}-\zeta^8\th^2\varphi v_{xxx}^2\]dxdt\\
\ns&\ds =\int_{Q} \[ -(\zeta^8\th^2\varphi)_x v_x(\cP v-v_{xxxx})+\frac{1}{2}\g\zeta^8(\th^2\varphi)_tv_x^2- (\zeta^8\th^2\varphi)_{x}v_{xx}v_{xxx}-\zeta^8\th^2\varphi v_{xxx}^2\]dxdt\\
\ns&\ds =\int_{Q} \[ -(\zeta^8\th^2\varphi)_x v_x\cP v-(\zeta^8\th^2\varphi)_{xx}v_xv_{xxx} -2(\zeta^8\th^2\varphi)_{x}v_{xx}v_{xxx}\\
\ns&\ds\q\q\q+\frac{1}{2}\g\zeta^8(\th^2\varphi)_tv_x^2-\zeta^8\th^2\varphi v_{xxx}^2\]dxdt.\\
\end{array}
\end{eqnarray*}
Then, noting that $\zeta\le 1$, we have
\begin{eqnarray}\label{c0201}
\begin{array}{rl}
&\ds  \int_{Q} \zeta^8 \l s^2\th^2\varphi  v_{xxx}^2dxdt\\
\ns&\ds \q\le C \int_{Q}\th^2|\cP v|^2dxdt +C\int_{Q}\th^2\zeta^4\(\zeta^2\l^3s^4\varphi^3v_{xx}^2+\l^5s^6\varphi^5v_{x}^2\)dxdt.
\end{array}
\end{eqnarray}
Multiplying $\g v_t+v_{xxxx}$ by $\zeta^6\th^2\varphi^3 v$ and integrating it over $Q$, we get
\begin{eqnarray*}
\begin{array}{rl}
&\ds  \int_{Q}  (\g v_t+v_{xxxx})\zeta^6\th^2\varphi^3 vdxdt\\
\ns&\ds =\int_{Q}  \[-\frac{1}{2}\g(\zeta^6\th^2\varphi^3)_t v^2-(\zeta^6\th^2\varphi^3)_{x}vv_{xxx}-\zeta^6\th^2\varphi^3v_xv_{xxx}\]dxdt\\
\ns&\ds =\int_{Q}  \[-\frac{1}{2}\g(\zeta^6\th^2\varphi^3)_t v^2-(\zeta^6\th^2\varphi^3)_{x}vv_{xxx}+(\zeta^6\th^2\varphi^3)_xv_xv_{xx}+\zeta^6\th^2\varphi^3v_{xx}^2\]dxdt.\\
\end{array}
\end{eqnarray*}
Then, for any $\e_1>0$, it holds that
\begin{eqnarray*}
\begin{array}{rl}
&\ds  \int_{Q}\zeta^6\l^3 s^4\th^2\varphi^3 v_{xx}^2dxdt\\
\ns&\ds \q\le C \int_{Q}\th^2|\cP v|^2dxdt +C\int_{Q}\th^2\zeta^2\(\l^7s^8\varphi^7v^2+\zeta^2\l^5s^6\varphi^5v_{x}^2\)dxdt\\
\ns&\ds \q\q+\e_1 \int_{Q} \zeta^8 \l s^2\th^2\varphi  v_{xxx}^2dxdt.
\end{array}
\end{eqnarray*}
Hence, for sufficiently small $\e_1$, by (\ref{c0201}),  we have
\begin{eqnarray}\label{c201}
\begin{array}{rl}
&\ds  \int_{Q}\zeta^6\l^3 s^4\th^2\varphi^3 v_{xx}^2dxdt\\
\ns&\ds \q\le C \int_{Q}\th^2|\cP v|^2dxdt +C\int_{Q}\th^2\zeta^2\(\l^7s^8\varphi^7v^2+\zeta^2\l^5s^6\varphi^5v_{x}^2\)dxdt.
\end{array}
\end{eqnarray}
Further, by (\ref{c0201}), it is easy to see that for any $\e_2>0$
\begin{eqnarray*}
\begin{array}{rl}
&\ds  \int_{Q}\zeta^4\l^5s^6\th^2\varphi^5v_{x}^2dxdt=\l^5s^6\int_{Q}\[-(\zeta^4\th^2\varphi^5)_xv_{x}v-\zeta^4\th^2\varphi^5vv_{xx}\]dxdt\\
\ns&\ds \q= \l^5s^6\int_{Q}\[\frac{1}{2}(\zeta^4\th^2\varphi^5)_{xx}v^2-\zeta^4\th^2\varphi^5vv_{xx}\]dxdt\\
\ns&\ds \q\le C\int_{Q}\zeta^2\l^7 s^8\th^2\varphi^7v^2dxdt+\e_2 \int_{Q}\zeta^6\l^3 s^4\th^2\varphi^3 v_{xx}^2dxdt\\
\end{array}
\end{eqnarray*}
Therefore, taking $\e_2$ small enough and combining (\ref{c201}), we end up with
\begin{eqnarray}\label{c00201}
\begin{array}{rl}
&\ds  \int_{Q}\zeta^4\l^5s^6\th^2\varphi^5v_{x}^2dxdt\\
\ns&\ds  \q\le C \int_{Q}\th^2|\cP v|^2dxdt +C\int_{Q}\zeta^2\l^7s^8\th^2\varphi^7v^2dxdt.
\end{array}
\end{eqnarray}
Combining (\ref{c0201}), (\ref{c201}) and (\ref{c00201}), we have
\begin{eqnarray}\label{c000201}
\begin{array}{rl}
&\ds  \int_{0}^T\int_{\omega_0}  \th^2\l s^2\varphi \( v_{xxx}^2+\l^2s^2\varphi^2 v_{xx}^2+\l^4s^4\varphi ^4v_x^2+\l^6s^6\varphi^6 v^2\)dxdt\\
\ns&\ds  \q\le C \int_{Q}\th^2|\cP v|^2dxdt+ C\int_{0}^T\int_{\omega}\l^7s^8\th^2\varphi^7v^2dxdt .
\end{array}
\end{eqnarray}
 Then,   combining  (\ref{c8}) and (\ref{c000201}), we can easily get (\ref{3a2-01}).

  \endpf

\section{Conditional stability in the inverse problem for the half-order time fractional Schr\"odinger equation }\label{section4}

As application of Theorem \ref{thm-0}, in this section, we will discuss an inverse source problem of a time fractional Schr\"odinger equation with half-order in time variable.
We consider the following time fractional Schr\"odinger equation:
   \bel{0a1-0}\left\{\ba{ll}\ds
(i\pa_t^{\frac{1}{2}}+\pa_x^2) z(t,x)=q(x)z+r(x)R(t,x),&\mbox{ in } Q,\\
\ns\ds z(0,x)=0,&\mbox{ in } \O,\\
\ns\ds z(t,0)=0,\q\pa_x z(t,0)=0,&\mbox{ on }(0,T).
\ea\right.\ee
In (\ref{0a1-0}),
 $q, R$ are given, $r$ is unknown and $\pa_t^{\frac{1}{2}}$ is the fractional order derivative in the Caputo sense which is defined by
 \begin{eqnarray*}\ba{ll}\ds
 \pa_t^\a z(t)=\frac{1}{\G(1-\a)}\int_0^t\frac{\pa_s z(s)}{(t-s)^\a}ds,\q 0<\a<1.
 \ea\end{eqnarray*}

 In the literature, there are numerous papers  dealing with the inverse problems for the Schr\"odinger equation with integer order time derivative, see \cite{BP, BC, Ben, E, HKSY, K, MOR} and the references therein. Recently, inverse problems for the fractional order heat equations draw more attention, there are  few interesting papers addressing this topic, not only in theoretical interests but also in applications (see \cite{CNYY, NSY} for example). In \cite{XCY, YZ}, the authors investigate two kinds of  conditional stability problems for  one dimensional heat equations with $1/2$-order time derivative. As far as we know, there is no references addressing inverse problem of the time fractional Schr\"odinger equation. The time fractional Schr\"odinger equation can be derived from the L\'evy-style quantum mechanics path (e.g. {\cite{MR2077514}}). We refer to \cite[ pp 20--23]{MR3380802}, and \cite{MR1809268} for the physical background of this equation.

 \subsection{Main Result}

Before giving the statement of our main result, we first  introduce some notations and assumptions.
We choose $\beta$ in (\ref{0a3-1}) satisfying
  \bel{0a4}
  x_0^2<\b\min\{t_0^2, (T-t_0)^2\}.
  \ee
By (\ref{0a4}) and the definition of $\psi$ in (\ref{0a3-1}), we know that if $\psi(t,x)>0$ for $0<x<L$, then $0<t<T$.

We set
 \bel{0a3x}\ba{ll}\ds
 $$
 \O_\e=Q_\e\cap\{t=t_0\},
 $$
  \ea\ee
where $Q_\e$ is given in  (\ref{0a3x-1}).
Then, it is easy to see that
\begin{eqnarray*}\left\{\ba{ll}\ds
 \O_\e=(0, x_0-\sqrt{\e}), &\mbox{ for } (x_0-L)^2<\e<x_0^2,\\
 \ns\ds \O_{3\e}=(0, x_0-\sqrt{3\e}), &\mbox{ for } (x_0-L)^2<\e<\frac{1}{3}x_0^2.
 \ea\right.\end{eqnarray*}

Next, by (\ref{0a4}), we know that for any $0<c<\e$, there exists a sufficiently small $\tau_0>0$  such that
\begin{eqnarray*}
\big\{t\ | 0\le x\le L, \ \psi(t,x)>c\big\}\subset(\tau_0, T-\tau_0).
\end{eqnarray*}
Therefore, $Q_{c}\cap\oO\subset (\tau_0, T-\tau_0)\t \oO$.

\bt\label{t01} Let $z\in C([0,T]; L^2(\O))\cap C^2((0,T); W^{2,\i}(\O)\cap H^4(\O))\cap C^3((0,T); L^2(\O))$ satisfy (\ref{0a1-0}) and $q, r\in W^{2,\i}(\O)$ satisfies $r(0)=\pa_x r(0)=0$.   Suppose that
 \begin{eqnarray*}\left\{\ba{ll}\ds
 R\in C^2((0,T); W^{2,\i}(\O))\cap C([0,T]; L^\i(\O)),\\
 \ns\ds \pa_t^{1/2}R\in C^2((0,T); L^\i(\O)),\q R(t_0,x)\neq 0,\q x\in\oO.
 \ea\right.\end{eqnarray*}
Then, for any $\e\in  ((x_0-L)^2, \frac{1}{3}x_0^2)$, there exist constant $C>0$ and $\kappa\in(0,1)$ such that
\begin{eqnarray*}
 ||r||_{H^2(\O_{3\e})}\le C||z(t_0,\cdot)||^\kappa_{H^4(\O_\e)}.
 \end{eqnarray*}
\et

\br
Though a suitable transformation can change the fractional operator $i\pa_t+\pa_x^2$ to $\pa_t+\pa_x^4$, a new term arising with singularity about the  time, (see (\ref{1102-3})), thus differ from inverse coefficient problems of equations with integer order time derivative, we need additional data $z(t_0,\cdot)$ to deal with inverse problem of fractional order equation.
\er

\ms

 \subsection{Proof of Theorem \ref{t01}}

 In this subsection, we will give the proof of Theorem \ref{t01}. To this aim, we first recall the following known results.

  \bl\label{thm-00} Let $a(x)\in C^1(\oO)$. Set $\Theta=e^{\l\psi(t_0,x)}$ with $\psi$ satisfying (\ref{0a3-1}) and $t_0\in(0, T)$. Then for any $u\in H_0^2(\O)$, there is a $\l_1>0$ such that for any $\l>\l_1$, the following inequality holds:
  \begin{eqnarray*}\ba{ll}\ds
  \int_{\O}\Th^2\( \l^3|u|^2+\l|\pa_x u|^2+\l^{-1} |\pa_x(a(x)\pa_x u)|^2\)dx\le C\int_{\O}\Th^2 |\pa_x(a(x)\pa_x u)|^2 dx.
  \ea\end{eqnarray*}
 \el

\ms

{\it Proof.} In \cite[Lemma 2.1]{FLZ},  by choosing  $m=1,\  \Psi=\Phi=0$,\  $\ell=\l|x-x_0|^2$, and $a^{11}=a(x)$,  a short calculation shows yields the desired result.\endpf

\bl\label{prop2} (\cite{XCY}) Let $z\in C[0,T]\cap W^{1,1}(0,T)$ and
\begin{eqnarray*}
z(0)=\pa_t^\b z(0)=0.
 \end{eqnarray*}
Then
 $$
\pa_t^{\a}\pa_t^{\b}z=\pa_t^{\a+\b}z,\q 0<\a+\b\le1.
$$
\el

\ms

{\bf Proof of Theorem \ref{t01}} We divide the proof into several steps.

\ms

{\bf Step 1: A transformation. } Set
\bel{1102-3}
\h z(t,x)=z(t,x)+\frac{2ir(x)R(0,x)t^{1/2}}{\G(1/2)}.
\ee
Then we can easily check that
\begin{eqnarray*}
\h z(0,x)=0
\end{eqnarray*}
and
\begin{eqnarray*}\ba{ll}\ds
\pa_t^{1/2}\h z(0,x)&\ds=\pa_t^{1/2}z(0,x)+\frac{2ir(x)R(0,x)}{\G(1/2)}\cdot\frac{\G(1/2)}{2} \\
\ns&\ds=\pa_t^{1/2}z(0,x)+i(i\pa_t^{1/2}+\pa_x^2-q(x))z(0,x) \\
\ns&\ds=\pa_t^{1/2}z(0,x)-\pa_t^{1/2}z(0,x)+i\pa_x^2z(0,x)-iq(x)z(0,x) \\
\ns&\ds=0.
\ea\end{eqnarray*}
Thus we obtain
\bel{1102-6}\ba{ll}\ds
\pa_t\h z&\ds=-i\pa_t^{1/2}(i\pa_t^{1/2}\h z) \\
\ns&\ds=-i\pa_t^{1/2}\[i\pa_t^{1/2}z+i\pa_t^{1/2}\(\frac{2irR(0,x)t^{1/2}}{\G(1/2)}\)\] \\
\ns&\ds=-i\pa_t^{1/2}[-\pa_x^2z+qz+rR-rR(0,x)] \\
\ns&\ds=(\pa_x^2-q)(i\pa_t^{1/2}z)-i\pa_t^{1/2}(rR) \\
\ns&\ds=(\pa_x^2-q)(-\pa_x^2z+qz+rR)-i\pa_t^{1/2}(rR) \\
\ns&\ds=-\pa_x^4z+(1+q)\pa_x^2z+2\pa_xq\pa_xz+(\pa_x^2q-q^2)z \\
\ns&\ds\q +(-i\pa_t^{1/2}+\pa_x^2-q)(rR).
\ea\ee
According to (\ref{1102-3}), we have
\bel{1102-7}\ba{ll}\ds
\pa_t\h z=\pa_tz+\frac{iR(0,x)}{\G(1/2)t^{1/2}}r.
\ea\ee
Combining (\ref{1102-6}) and (\ref{1102-7}), we obtain
\begin{eqnarray*}\ba{ll}\ds
(\pa_t+\pa_x^4)z&\ds=(1+q)\pa_x^2z+2\pa_xq\pa_xz+(\pa_x^2q-q^2)z \\
\ns&\ds\q +\pa_x(R\pa_xr)+\pa_xR\pa_xr+(-i\pa_t^{1/2}R+\pa_x^2R-qR-\frac{iR(0,x)}{\G(1/2)t^{1/2}})r.
\ea\end{eqnarray*}

Thus,
  \bel{f04-06}\ba{ll}\ds
 \cP z= \pa_tz+\pa_x^4z=\h{\cP}z+L r,
  \ea\ee
where
  \bel{f04-06x1}\left\{\ba{ll}\ds
\h{\cP}z\=(1+q)\pa_x^2z+2\pa_xq\pa_xz+(\pa_x^2q-q^2)z,\\
\ns\ds
Lr\=\pa_x(R\pa_xr)+\pa_xR\pa_xr+(-i\pa_t^{1/2}R+\pa_x^2R-qR-\frac{iR(0,x)}{\G(1/2)t^{1/2}})r.
  \ea\right.\ee

For simplicity, we set
\bel{f04-06x4}\ba{ll}\ds
\(\cP_0z\)(x)\=\(\cP z\)(t_0,x),\q \(\h{\cP}_0z\)(x)\=\(\h{\cP}z\)(t_0,x),\q \(L_0r\)(x)\=\(Lr\)(t_0,x).
\ea\ee

  Taking $t=t_0$ in equation (\ref{f04-06}) we have
 \begin{eqnarray*}\ba{ll}\ds
\cP_0 z=\h{\cP}_0z+L_0r.
   \ea\end{eqnarray*}

\ms

{\bf Step 2: Estimation of $r$ in $H^2$-norm localized in $\O_{2\e}$. }

Let $\wt{\chi}\in C^\i(\cl{\O}; [0,1])$ satisfying
 \begin{eqnarray*}
  \wt{\chi}=\left\{\ba{ll}\ds 1&\mbox{ in } \O_{2\e},\\
  \ns\ds
  0 &\mbox{ in } \O\setminus\cl{\O}_{\e},\ea\right.
   \end{eqnarray*}
where $\O_{\e}\=Q_\e\cap\{t=t_0\}$.

Next we define $\wt{r}=r\wt{\chi}$, by (\ref{f04-06}), (\ref{f04-06x1}) and (\ref{f04-06x4}), we have that
\begin{eqnarray*}\ba{ll}\ds
  (L_0\wt{r})(x)=\wt{\chi}L_0r(x)+\wt{g},
   \ea\end{eqnarray*}
   where
  $$
   \wt{g}=2R(t_0,x)(\pa_x\wt{\chi})\pa_xr(x)+R(t_0,x)(\pa_x^2\wt{\chi})r(x)+2\pa_xR(t_0,x)(\pa_x\wt{\chi})r(x).
 $$

  In Lemma \ref{thm-00}, we choose $a(x)=R(t_0,x),\  u=\wt r$, then there is a $\l_1>0$ and $C>0$, for any $\l>\l_1$, we have that
 \begin{eqnarray*}\ba{ll}\ds
    \int_\O(\l|\pa_x \wt r|^2+\l^3|\wt r|^2+\l^{-1}|\pa_{xx}\wt r|^2)e^{2\l\psi(t_0,x)}dx\\
    \ns\ds\le C\int_\O\|\pa_{x}(R(t_0,x)\pa_{x} \wt r)\|^2e^{2\l\psi(t_0,x)}dx\\
    \ns\ds\le C\int_\O\(|\wt{\chi}L_0r(x)|^2+|\pa_{x}\wt r|^2+|\wt r|^2+|\wt g|^2\)e^{2\l\psi(t_0,x)}dx.
   \ea\end{eqnarray*}
 Therefore, there exists a $\l_2>\l_1$, for any $\l>\l_2$, the following inequality holds:
   \bel{1028-07}\ba{ll}\ds
    \int_\O(\l|\pa_{x}\wt r|^2+\l^3|\wt r|^2+\l^{-1}|\pa_{xx}\wt r|^2)e^{2\l\psi(t_0,x)}dx&\ds\le  C\int_\O\(|\wt{\chi}L_0r(x)|^2+|\wt g|^2\)e^{2\l\psi(t_0,x)}dx.
   \ea\ee
Remember that  $\wt r=r\wt \chi$ and $L_0 r=\cP_0z-\h{\cP}_0z$, and noting $|\wt g|\neq 0$ only if $\psi\leq 2\e$, by (\ref{1028-07}), we have
  \bel{f04-25}\ba{ll}\ds
  \q\int_{\O_{2\e}}(\l|\pa_x r|^2+\l^3|r|^2+\l^{-1}|r_{xx}|^2)e^{2\l\psi(t_0,x)}dx \\
\ns\ds\le C\int_{\O_\e}|\wt{\chi}\cP_0z-\wt{\chi}\h{\cP}_0z|^2e^{2\l\psi(t_0,x)}dx+C\int_{\O_\e\setminus\O_{2\e}}|\wt g|^2e^{2\l\psi(t_0,x)}dx\\
\ns\ds\le C\int_{\O_{\e}}(|\pa_tz(t_0,x)|^2+\sum_{j=0}^4|\pa_x^jz(t_0,x)|^2)e^{2\l\psi(t_0,x)}dx+Ce^{4\l\e}.
  \ea\ee

{\bf Step 3: Estimation of $\ds\int_{\O_{\e}}|\pa_tz(t_0,x)|^2e^{2\l\psi(t_0,x)}dx$. }

Note that $Q_{\e}\cap\{t=0\}$ is empty by assumption (\ref{0a4}). Let $\cl\chi(t,x)\in C^\i(\cl{Q}; [0,1])$ satisfying
  \bel{f04-100}\ba{ll}\ds
\cl\chi=\begin{cases} 1\q \mbox{ in }\q Q_{\e},\\0\q \mbox{ in }\q Q\setminus\cl{Q}_{\frac{\e}{2}},\end{cases}
   \ea\ee
where $Q_{\frac{\e}{2}}\=\{(t,x)\in Q;\  \psi(t,x)>\frac{\e}{2}\}$.
\ms

By (\ref{f04-100}), we have $\cl\chi(t_0,x)=1$ for $x\in\O_{\e}$ and $\cl\chi=0$ on $\pa Q_{\frac{\e}{2}}\cap\{x>0\}$. Therefore
 \begin{eqnarray*}\ba{ll}&\ds
  \int_{\O_{\e}}|\pa_tz(t_0,x)|^2e^{2\l\psi(t_0,x)}dx\\
  \ns&\ds\le\int_{\O_{\frac{\e}{2}}}|\cl\chi\pa_tz(t_0,x)|^2e^{2\l\psi(t_0,x)}dx=\int_{Q_{\frac{\e}{2}}\cap\{t<t_0\}}\pa_t(|\cl\chi\pa_tz|^2e^{2\l\psi})dxdt \\
  \ns&\ds\le C\int_{Q_{\frac{\e}{2}}}(\cl\chi|\pa_t\chi||\pa_tz|^2+\cl\chi^2|\pa_tz||\pa_t^2z|+\l\cl\chi^2|\pa_tz|^2)e^{2\l\psi(t,x)}dxdt \\
  \ns&\ds\le C\int_{Q_{\frac{\e}{2}}}\l(|\pa_tz|^2+|\pa_t^2z|^2)e^{2\l\psi(t,x)}dxdt.
   \ea\end{eqnarray*}
   
   Let $\chi(t,x)\in C^\i(\cl{Q}; [0,1])$ satisfying
  \bel{f04-10}\ba{ll}\ds
\chi=\begin{cases} 1\q \mbox{ in }\q Q_{\frac{\e}{2}},\\0\q \mbox{ in }\q Q\setminus\cl{Q}_{\frac{\e}{3}}.\end{cases}
   \ea\ee
   
By (\ref{f04-06})  and (\ref{f04-06x1}), we know that
  \begin{eqnarray*}\ba{ll}\ds
\cP(\chi\pa_tz)=\chi\pa_t(Lr)+\h{\cP}(\chi\pa_tz)+\wt{g_1},\q
\cP(\chi\pa_t^2z)=\chi\pa_t^2(Lr)+\h{\cP}(\chi\pa_t^2z)+\wt{g_2},
   \ea\end{eqnarray*}
where $\wt{g_j}$ are linear combinations of $\pa_x^i\pa_t^jz, i=0,1,2,3, j=1,2$ whose coefficients contain $\pa_t\chi,\pa_x^i\chi,i=0, 1,2,3,4$ as multiplicative factors. By the regularity of $z$ and the definition of $\chi$, we note that $\wt{g_j}\in L^2(Q_{\e})$.

In Theorem \ref{thm-0}, we choose $v=\chi\pa_tz$. Then there is a $\l_0>0$, for any $\l>\l_0$, we have that
  \bel{f04-14}\ba{ll}\ds
  \q\int_{Q_{\frac{\e}{3}}}\( \l^7|\chi\pa_tz|^2+\l^5|\pa_x(\chi\pa_tz)|^2+\l^3|\pa_x^2(\chi\pa_tz)|^2+\l|\pa_x^3(\chi\pa_tz)|^2\)e^{2\l\psi}dxdt\\
  \ns\ds\le C\int_{Q_{\frac{\e}{3}}}|\cP(\chi\pa_tz)|^2e^{2\l\psi}dxdt\\
  \ns\ds\le C\int_{Q_{\frac{\e}{3}}}\(|\chi\pa_t(Lr)+\sum_{j=0}^2|\pa_x^j(\chi \pa_t z)|^2\)e^{2\l\psi}dxdt+C\int_{Q_{\frac{\e}{3}}\setminus Q_{\frac{\e}{2}}}|\wt{g_1}|^2e^{2\l\psi}dxdt.
   \ea\ee
Therefore, there exists a $\l_3>\l_0$, for any $\l>\l_3$,
  \bel{f04-14y}\ba{ll}\ds
  \q\int_{Q_{\frac{\e}{3}}}\l^7|\chi\pa_t z|^2 e^{2\l\psi}dxdt
\le C\int_{Q_{\frac{\e}{3}}}|\pa_t(Lr)|^2e^{2\l\psi}dxdt+C\int_{Q_{\frac{\e}{3}}\setminus Q_{\frac{\e}{2}}}|\wt{g_1}|^2e^{2\l\psi}dxdt.
   \ea\ee
   Similarly, in Theorem \ref{thm-0}, we choose $v=\chi\pa^2_tz$, then there exists a $\l_4>\l_0$, for any $\l>\l_4$,
     \bel{f04-14y1}\ba{ll}\ds
  \q\int_{Q_{\frac{\e}{3}}}\l^7|\chi\pa^2_t z|^2 e^{2\l\psi}dxdt
\le C\int_{Q_{\frac{\e}{3}}}|\pa^2_t(Lr)|^2e^{2\l\psi}dxdt+C\int_{Q_{\frac{\e}{3}}\setminus Q_{\frac{\e}{2}}}|\wt{g_2}|^2e^{2\l\psi}dxdt.
   \ea\ee
Noting that $r(x)$ independent of $t$, then we have
  \bel{f04-16}\ba{ll}\ds
|\pa_t(Lr)|^2+|\pa^2_t(Lr)|^2\le C(|\pa_x^2r|^2+|\pa_xr|^2+|r|^2).
\ea\ee

Further,
  \bel{21121001}\ba{ll}\ds
\int_{Q_{\frac{\e}{2}}}\l(|\pa_tz|^2+|\pa_t^2z|^2)e^{2\l\psi}dxdt \le \int_{Q_{\frac{\e}{3}}}\l(|\chi\pa_tz|^2+|\chi\pa_t^2z|^2)e^{2\l\psi}dxdt.
\ea\ee

Combining (\ref{f04-14y})--(\ref{21121001}), taking $\l_5=\max\{\l_3, \l_4\}$, for any $\l>\l_5$, we get
  \bel{f04-17}\ba{ll}\ds
  \int_{Q_{\frac{\e}{2}}}\l(|\pa_tz|^2+|\pa_t^2z|^2)e^{2\l\psi}dxdt\le \frac{C}{\l^6}\int_{Q_{\frac{\e}{3}}}(|\pa_x^2r|^2+|\pa_xr|^2+|r|^2)e^{2\l\psi}dxdt+Ce^{\l\e}.
   \ea\ee

Finally, by (\ref{f04-25}) and (\ref{f04-17}), we conclude that there is a $\l_6=\max\{\l_2, \l_5\}$, for any $\l>\l_6$, we have
  \bel{f04-27}\ba{ll}\ds
  \int_{\O_{2\e}}(\frac{1}{\l}|\pa_x^2r|^2+\l|\pa_x r|^2+\l^3|r|^2)e^{2\l\psi(t_0,x)}dx \\
  \ns\ds\le C\int_{\O_\e}\sum_{j=0}^4|\pa_x^jz(t_0,x)|^2e^{2\l\psi(t_0,x)}dx+Ce^{4\l\e}\\
  \ns\ds\q +\frac{C}{\l^6}\int_{Q_{\frac{\e}{3}}}(|\pa_x^2r|^2+|\pa_xr|^2+|r|^2)e^{2\l\psi}dxdt.
\ea\ee

\ms

{\bf Step 4: End of the proof. }

 On the one hand, in $\cl{Q_\e}$, it is easy to see that $\e\le\psi(t,x)\le|x_0|^2$.  Therefore,
  \bel{f04-28}\ba{ll}\ds
\int_{\O_\e}\sum_{j=0}^4|\pa_x^jz(t_0,x)|^2e^{2\l\psi(t_0,x)}dx \le Ce^{2\l|x_0|^2}\int_{\O_\e}\sum_{j=0}^4|\pa_x^jz(t_0,x)|^2dx.
\ea\ee
On the other hand,
  \bel{f04-29}\ba{ll}\ds
\frac{C}{\l^6}\int_{Q_{\frac{\e}{3}}}(|\pa_x^2r|^2+|\pa_x r|^2+|r|^2)e^{2\l\psi(t,x)}dxdt \\
\ns\ds\le \frac{C}{\l^6}\(\int_{\O_{\frac{\e}{3}}\setminus\O_{2\e}}+\int_{\O_{2\e}}\)(|\pa_x^2r|^2+|\pa_x r|^2+|r|^2)e^{2\l\psi(t_0,x)}dx.
\ea\ee
Since $r\in W^{2,\i}(\O)$ and $\psi(t_0,x)\le2\e$ in $\O_{\frac{\e}{3}}\setminus\cl{\O_{2\e}}$ we can estimate that
  \bel{f04-30}\ba{ll}\ds
\int_{\O_{\frac{\e}{3}}\setminus\O_{2\e}}(|\pa_x^2r|^2+|\pa_x r|^2+|r|^2)e^{2\l\psi(t_0,x)}dx\le Ce^{4\l\e}.
\ea\ee
Therefore by (\ref{f04-27})--(\ref{f04-30}), we conclude that there is a $\l_7>\l_6$, for any $\l>\l_7$, we get
  \bel{f04-31}\ba{ll}\ds
\int_{\O_{2\e}}(\frac{1}{\l}|\pa_x^2r|^2+\l|\pa_x r|^2+\l^3|r|^2)e^{2\l\psi(t_0,x)}dx\le Ce^{4\l\e}+Ce^{2\l|x_0|^2}||z(t_0,\cdot)||_{H^4(\O_\e)} .
\ea\ee
Further, since $\O_{3\e}\subset\O_{2\e}$ and $\psi(t_0,x)>3\e$ in $\O_{3\e}$, then
  \bel{f04-32}\ba{ll}\ds
  \int_{\O_{2\e}}(\frac{1}{\l}|\pa_x^2r|^2+\l|\pa_x r|^2+\l^3|r|^2)e^{2\l\psi(t_0,x)}dx \\
  \ns\ds\ge\int_{\O_{3\e}}(\frac{1}{\l}|\pa_x^2r|^2+\l|\pa_x r|^2+\l^3|r|^2)e^{2\l\psi(t_0,x)}dx \\
  \ns\ds\ge e^{6\l\e}\int_{\O_{3\e}}(\frac{1}{\l}|\pa_x^2r|^2+\l|\pa_xr|^2+\l^3|r|^2)dx.
  \ea\ee
Combining (\ref{f04-31}) and (\ref{f04-32}), we have
  \bel{f04-33}\ba{ll}\ds
\int_{\O_{3\e}}(\frac{1}{\l}|\pa_x^2r|^2+\l|\pa_x r|^2+\l^3|r|^2)dx\le Ce^{-2\l \tau}+Ce^{C\l}||z(t_0,\cdot)||_{H^4(\O_\e)}^2,
\ea\ee
where $\tau={3\e}-{2\e}>0$. Minimizing the right-hand side of (\ref{f04-33}), by choosing $\l$ large enough, we get
 \begin{eqnarray*}\ba{ll}\ds
||r||_{H^2(\O_{3\e})}\le C||z(t_0,\cdot)||_{H^4(\O_{\e})}^\k.
\ea\end{eqnarray*}
This completes the proof of Theorem \ref{t01}. \endpf


\begin{thebibliography}{99}

\bibitem{BP} L. Baudouin, J.-P. Puel, \it Uniqueness and stability in an inverse problem for the Schr\"odinger equation, \sl
Inverse Probl., \rm {\bf 18}(2002), 1537--1554.

\bibitem{BC} M.~Bellassoued, M.~Choulli, \it Stability estimate for an inverse problem for the magnetic Schr\"odinger
equation from the Dirichlet-to-Neumann map, \sl J. Funct. Anal., \rm {\bf 91}(2010), 161--195.

\bibitem{Ben}H.~Ben Joud, \it A stability estimate for an inverse problem for the Schr\"odinger equation in a magnetic
field from partial boundary measurements, \sl Inverse Probl., \rm {\bf 25}(2009), 45012--45034.



\bibitem{c} T. Carleman, {\it Sur un probl\`{e}me d'unicit\'{e} pur les syst\`{e}mes
d'\'{e}quations aux d\'{e}riv\'{e}es partielles \`{a} deux variables ind\'{e}pendantes},
Ark. Mat. Astr. Fys., \rm{\bf 26}(1939), 1--9.


\bibitem{CM} E.~Cerpa, A.~Mercado, \it Local exact controllability to the trajectories of the 1-D Kuramoto-Sivashinsky equation, \sl J. Differential Equations, \rm {\bf 250}(2011), 2024--2044.

\bibitem{CNYY} J.~Cheng, J.~Nakagawa, M.~Yamamoto and T.~Yamazaki, \it Uniqueness in an inverse problem for a one-dimensional fractional diffusion equation, \sl Inverse Probl., \rm {\bf 25}(2009), 115002.



\bibitem{E} G.~Eskin, \it Inverse boundary value problems and the Aharonov-Bohm effect, \sl Inverse Probl., \rm {\bf 19}(2003), 49--62.

\bibitem{FLZ}X. Fu, Q. L\"{u} and X. Zhang, \it Carleman estimates for second order partial differential operators an applications, a unified approach, \sl \upshape Springer, \rm 2019.

    \bibitem{FI} A.~V.~Fursikov and O.~Yu.~Imanuvilov, \sl
Controllability of Evolution Equations, \rm Lecture Notes Series
{\bf 34}, Research Institute of Mathematics, Seoul National
University, Seoul, Korea, \rm 1996.

\bibitem{GK} S.~Guerrero and K.~Kassab, \it Carleman estimate and null controllability of a fourth order
parabolic equation in dimension $n\ge 2$,  \sl J. Math. Pures Appl.,  \rm {\bf 121}(2019), 135--161.

\bibitem{MR3380802}B.~Guo, X.~Pu, and F.~Huang, \it Fractional partial differential
	equations and their numerical solutions, \sl World Scientific Publishing Co.
Pte. Ltd., Hackensack, NJ, \rm 2015. \newblock Originally published by Science Press in 2011.

 \bibitem{HK} X.~Huang, A.~Kawamoto, \it Inverse problems for a half-order time-fractional diffusion equation in arbitrary dimension by Carleman estimates, \rm https://arxiv.org/abs/2010.10082.

      \bibitem{HKSY} X.~Huang, Y.~Kian, \'E.~Soccorsi and M.~Yamamoto, \it Carleman estimate for the Schr\"odinger equation and application to magnetic inverse problems,  \sl J. Math. Anal. Appl.,  \rm {\bf 474}(2019), 116--142.

     \bibitem{Hormander} L. H\"ormander, \it  Uniqueness theorems for second
order elliptic differential equations. \sl Commun.
Part. Diff. Eq., \rm \textbf{8}(1983), 21--64.

\bibitem{Isakov2} V. Isakov, \it Inverse Problems for Partial Differential Equations, \sl Springer,  \rm 2006.

\bibitem{Kn}  C. E. Kenig, \it  Carleman estimates, uniform Sobolev inequalities for second order differential operators, and unique continuation theorems, \sl Proceedings of the International Congress of Mathematicians, \rm {\bf 2}(1986), 948--960.

    \bibitem{Kli} M.~V.~Klibanov, \it Carleman estimates for global
uniqueness, stability and numerical methods for
coefficient inverse problems. \sl J. Inverse
Ill-Posed Probl. \rm {\bf 21}(2013), 477--560.


 \bibitem{K}Y.~Kian, Q.~Phan, E.~Soccorsi, \it H\"older stable determination of a quantum scalar potential in
unbounded cylindrical domains, \sl J. Math. Anal. Appl.,  \rm {\bf 426}(2015), 194--210.



  \bibitem{MOR} A.~Mercado, A.~Osses, L.~Rosier, \it Inverse problems for the Schr\"odinger equation via Carleman inequalities with degenerate weights, \sl Inverse Probl., \rm {\bf 24}(2008), 015017.

\bibitem{MR1809268}
 R.~Metzler and J.~Klafter, \it The random walk's guide to anomalous
	diffusion: a fractional dynamics approach, \sl Phys. Rep., \rm {\bf 339}(2000), 1--77.

\bibitem{MR2077514} M.~Naber, \it Time fractional {S}chr\"{o}dinger equation, \sl J. Math. Phys., \rm {\bf45}(2004), 3339--3352.

\bibitem{NSY} J.~Nakagawa, K.~Sakamoto and M.~Yamamoto, \it Overview to mathematical analysis for fractional diffusion equations--new mathematical aspects motivated by industrial collaboration, \sl J. Math-for-Industry A, \rm {\bf 2}(2010 A-10), 98--108.


\bibitem{XCY} X. Xu, J. Cheng and M. Yamamoto,  \it Carleman estimate for a fractional diffusion equation with half order
and application, \sl Appl. Anal., \rm {\bf 90}(2011) 1355--1371.

\bibitem{Y} M.~Yamamoto, \it Carleman estimates for parabolic equations and applications, \sl Inverse Probl., \rm {\bf 25}(2019), 123013.

\bibitem{YZ} M. Yamamoto and Y. Zhang,  \it Conditional stability in determining a zeroth-order coefficient in a half-order fractional diffusion equation by a Carleman estimate, \sl Inverse Probl., \rm {\bf 28}(2012), 105010.



\bibitem{Zheng} C.~Zheng, \it Inverse problem for the fourth order Schr\"odinger equation on a finite domain, \sl  Math. Control Relat. Fields, \rm {\bf 5}(2015), 177--189.

\bibitem{ZZ}Z.~Zhou, \it Observability estimate and null controllability for one-dimensional fourth order parabolic equation, \sl Taiwan. J.
Math., \rm {\bf16}(2012), 1991--2017.

\end{thebibliography}
\end{document}